\documentclass[11pt,letterpaper, reqno]{amsart}
\usepackage{mathrsfs}
\usepackage[active]{srcltx}
\usepackage{mathrsfs}
\usepackage{longtable}
\usepackage{cancel}
\usepackage[normalem]{ulem}

\usepackage{xcolor}
\definecolor{bluecite}{HTML}{0875b7}
\usepackage[unicode=true,
bookmarksopen={true},
pdffitwindow=true,
colorlinks=true,
linkcolor=bluecite,
citecolor=bluecite,
urlcolor=bluecite,
hyperfootnotes=false,
pdfstartview={FitH},
pdfpagemode= UseNone]{hyperref}

\newcommand{\ds}{\displaystyle}

\newtheorem{proposition}{Proposition}[section]
\newtheorem{theorem}{Theorem}[section]
\newtheorem{lemma}{Lemma}[section]
\newtheorem{corollary}{Corollary}[section]

\newtheorem{remark}{Remark}[section]
\numberwithin{equation}{section}

\usepackage{bbm}  
\usepackage{dsfont}

\begin{document}
	\title[Logarithmic-Sobolev inequalities on submanifolds]
	%{Asymptotic volume ratios characterize sharp constants in Sobolev inequalities on Riemannian manifolds with nonnegative Ricci curvature 
		%}
	%Explicit sharp constants in geometric and functional inequalities on Riemannian manifolds with nonnegative Ricci curvature
	%{Isoperimetric inequalities in spaces with nonnegative Ricci curvature and volume growth at infinity
		%}
	{Logarithmic-Sobolev inequalities on non-compact Euclidean submanifolds: sharpness and rigidity 
		%
	%	Sharp logarithmic-Sobolev inequalities on complete %Euclidean submanifolds
	}

	\author{Zolt\'an M. Balogh and Alexandru Krist\'aly}	
	
%	\thanks{Z. M. Balogh was
%		supported by the Swiss National Science Foundation, Grant Nr. {200020\_191978}.  A. Krist\'aly  was supported by the UEFISCDI/CNCS grant PN-III-P4-ID-PCE2020-1001.}

	\address{Mathematisches Institute,
		Universit\"at Bern,
		Sidlerstrasse 5,
		3012 Bern, Switzerland}
	
	\email{zoltan.balogh@math.unibe.ch}
	
	\address{Department of Economics, Babe\c s-Bolyai University, Str. Teodor Mihali 58-60, 400591 Cluj-Napoca,
		Romania \& Institute of Applied Mathematics, \'Obuda
		University, B\'ecsi \'ut 96/B, 1034
		Budapest, Hungary}
	
	%\curraddr{Department of Mathematics and Statistics,
		%Case Western Reserve University, Cleveland, Ohio 43403}
	\email{alexandru.kristaly@ubbcluj.ro; kristaly.alexandru@uni-obuda.hu}

	
	
	\thanks{Z. M. Balogh is supported by the Swiss National Science Foundation, Grant nr.\ 200020-191978 and Grant nr.\ 200021-228012.\  A.\ Krist\'aly is  supported by the
		Excellence Researcher Program \'OE-KP-2-2022 of \'Obuda University, Hungary.}
	
	\begin{abstract}
		{\footnotesize \noindent 
		The paper is devoted to provide  Michael--Simon-type $L^p$-logarith\-mic-Sobolev inequalities on  complete,  not necessarily compact
		$n$-dimensional submanifolds $\Sigma$ of the Euclidean space $\mathbb R^{n+m}$.  Our first result, stated for $p=2$, is sharp, it is valid on general submanifolds, and it involves the mean curvature of $\Sigma$.  It  implies in particular the main result of S.\ Brendle [\textit{Comm.\ Pure Appl.\ Math.}, 2022]. In addition, it turns out that equality can only occur if and only if $\Sigma$ is isometric to the Euclidean space
		 $\mathbb R^{n}$ and the extremizer is a Gaussian. The second result is a general $L^p$-logarithmic-Sobolev inequality for $p\geq 2$ on Euclidean submanifolds with constants that are codimension-free in case of minimal submanifolds. In order to prove the above results -- especially, to deal with the equality cases --  we elaborate the theory of optimal mass transport on submanifolds between measures that are not necessarily compactly supported.
	Applications are provided to sharp hypercontractivity estimates of Hopf--Lax semigroups on  submanifolds. The first hypercontractivity estimate is for general submanifolds with bounded mean curvature vector, the second one is for self-similar shrinkers endowed with the natural Gaussian measure. The equality cases are characterized here as well.      
		}
	\end{abstract}

%By using optimal transport theory, we prove Michael--Simon-type $L^p$-logarith\-mic-Sobolev inequalities on %non-compact 
%$n$-dimensional submanifolds $\Sigma$ of the Euclidean space $\mathbb R^{n+m}$.  Our first result, stated for $p=2$, is sharp, it is valid on general submanifolds, and it involves the mean curvature of $\Sigma$.  This  implies in particular the main result of S.\ Brendle [\textit{Comm.\ Pure Appl.\ Math.}, 2023]. Considering the equality case in this result, we find that equality can only occur iff $\Sigma$ is isometric to the Euclidean space
%$\mathbb R^{n}$ and the \textcolor{red}{extremizer} is a Gaussian. The second result is a general $L^p$-logarithmic-Sobolev inequality for $p\geq 2$ on Euclidean submanifolds with constants that are codimension-free in case of minimal submanifolds. 
%%		We also characterize the equality cases for $p=2$, which occur if and only of the submanifolds are isometric to an Euclidean space and the extremals are Gaussian bubbles. 
%Applications are provided to sharp hypercontractivity estimates of Hopf--Lax semigroups on submanifolds. The first estimate is for general submanifolds with bounded mean curvature vector, the second one is for self-similar shrinkers endowed with the natural Gaussian measure. \textcolor{red}{The equality cases are also characterized.}

	%\vspace*{-1cm} 
	\maketitle
	
	%\tableofcontents

%	\vspace*{-0.5cm}
	\section{Introduction}

The isoperimetric inequality (equivalently, the $L^1$-Sobolev inequality) on Euclidean submanifolds has been established by Michael and Simon \cite{MS} and independently by Allard \cite{Allard}; this inequality is  codimen\-sion-free,  containing the mean curvature of the submanifold, but is not sharp. By a standard argument, this result implies an $L^p$-Sobolev inequality valid for any $p\in [1,n)$,  where $n$ is the dimension of the submanifold,   see e.g.\ Cabr\'e and Miraglio \cite{Cabre-Miraglio}. For a proof of the Michael--Simon--Sobolev inequality based on optimal  transport,  see Castillon \cite{C}. 

Very recently, Brendle \cite{Brendle-1} and also Brendle and Eichmair \cite{Brendle-Eichmair} proved a Michael--Simon-type isoperimetric inequality, which turns out to be sharp for every Euclidean submanifold of codimension two. This result, when combined with the classical rearrangement technique of Aubin \cite{Aubin} and Talenti \cite{Talenti} and the co-area formula, implies a sharp P\'olya--Szeg\H o inequality in the setting of  {\it minimal submanifolds}  of
codimension at most two. This in turn implies the sharp $L^p$-Sobolev inequality for every $p\in [1,n)$, see Brendle \cite[Theorem 5.8]{{Brendle-Toulouse}}. In the same geometric setting further sharp Sobolev inequalities (logarithmic-Sobolev, Gagliardo--Nirenberg--Sobolev) can be established.
% via the sharp P\'olya--Szeg\"o inequality. 
 On the other hand, it seems that a similar approach cannot be used to establish sharp and codimension-free Sobolev inequalities on general (non-minimal) Euclidean submanifolds.  

A prominent subclass of Sobolev inequalities -- not necessarily in the submanifold setting --  is represented by logarithmic-Sobolev inequalities. Motivation to study such inequalities arise from various applications: 
	monotonicity formula for the entropy under Ricci and heat flows, see Perelman \cite{Perelman} and Ni \cite{Ni}; monotonicity formula for mean curvature flow, see Huisken \cite{Huisken}; equilibrium for spin systems, see Guionnet and Zegarlinski \cite{GZ}; 
	quantum field theory, see  Glimm and Jaffe \cite{GJ}; hypercontractivity estimates for Hopf--Lax semigroups, see e.g.\ Bobkov,  Gentil and Ledoux \cite{BobkovGL}, Gentil \cite{Gentil}, Otto and Villani \cite{OV1, OV2}, Balogh, Krist\'aly and Tripaldi \cite{BKT}.  
	
		The  well-known form of the  Euclidean $L^p$-logarithmic-Sobolev inequality  in $\mathbb R^n$ for $1 < p < n$ 
	has been established by  Del Pino and Dolbeault \cite{delPinoDolbeault-2},  and first by Weissler \cite{Weissler} for $p=2$. It states that  every function  $f \in {W}^{1,p}(\mathbb R^n)$ with $\displaystyle\int_{\mathbb R^n} |f|^p = 1$ verifies
	\begin{equation}\label{e-sharp-log-Sobolev}
		\int_{\mathbb R^n}  |f|^p\log |f|^p \leq \frac{n}{p}	\log\left(\mathcal
		L_{p,n}\displaystyle\int_{\mathbb R^n} |\nabla f|^p \right)
	\end{equation}
	with the sharp  constant 
	$$		\mathcal
	L_{p,n}={\small
		\frac{p}{n}\left(\frac{p-1}{e}\right)^{p-1}\left(\omega_n{\Gamma\left(\frac{n}{q}+1\right)}\right)^{-\frac{p}{n}}},   
	$$
	%	\end{equation}
 where $\Gamma$ is the usual  Gamma-function, $q=\frac{p}{p-1}$ is the conjugate of $p$, and $\omega_n=\frac{\pi^\frac{n}{2}}{\Gamma(\frac{n}{2}+1)}$ is the volume of the unit ball in $\mathbb R^n$. 
The constant
$\mathcal L_{p,n}$ is sharp in \eqref{e-sharp-log-Sobolev} and 
%by density arguments, \eqref{e-sharp-log-Sobolev} can be extended to the space ${W}^{1,p}(\mathbb R^n)$; moreover, 
equality holds in \eqref{e-sharp-log-Sobolev} (see e.g.\ in \cite{BDK}) if and only if $f\in {W}^{1,p}(\mathbb R^n)$ belongs -- up to translations -- to the class of Gaussians of the form $f_\lambda(x)=\lambda^{\frac{n}{pq}}\left(\omega_n\Gamma\left(\frac{n}{q}+1\right)\right)^{-1/p}e^{-\frac{\lambda}{p} |x|^q}$, $x\in \mathbb R^n$ and $\lambda>0.$ 
	
 As might be expected, the situation is more delicate in the submanifold setting. The first main contribution is due to Ecker \cite{Ecker}, who proved a codimension-free, but non-sharp $L^2$-logarithmic-Sobolev inequality on Euclidean submanifolds. Very recently,  by using the ABP-method, Brendle \cite{Brendle-log} provided the sharp, codimension-free version of Ecker's $L^2$-logarithmic-Sobolev inequality. This result states that if $\Sigma$ is a \textit{compact} $n$-dimensional submanifold of $\mathbb R^{n+m}$ without boundary and $g$ is any positive smooth function on $\Sigma$, then 
	\begin{equation}\label{Brendle-initial-inequality}
		\int_\Sigma g\left(\log g + n + \frac{n}{2}\log(4\pi)\right)-\int_\Sigma\frac{|\nabla^\Sigma g|^2}{g} -\int_\Sigma |H|^2g\leq \int_\Sigma g \log\left(\int_\Sigma g\right).
	\end{equation}
	Here, $H$ stands for the mean curvature vector of $\Sigma$,  $\nabla^\Sigma$ is the gradient associated to $\Sigma$, and the integrals are  considered with respect to the natural canonical measure $d{\rm vol}_\Sigma$ on $\Sigma$.  
	It is remarkable that in contrast to Brendle's Michael--Simon-type isoperimetric inequality (that is up until now sharp only for $m= 2$), inequality \eqref{Brendle-initial-inequality} is sharp for {\it any} codimension $m\geq 1$.
	
	The main purpose of the present paper is to use optimal transport theory to prove  counterparts of the logarithmic-Sobolev inequality \eqref{e-sharp-log-Sobolev}  on \textit{complete}, not necessarily compact $n$-dimensional Euclidean submanifolds of $\mathbb R^{n+m}$. This problem has been recently considered by Wang, Xia and Zhang  \cite{Wang-Xia-Zhang} for the special case of \textit{non-compact self-shrinkers} by the so-called ABP method used also by Brendle.  Unlike in the compact case, see Brendle \cite{Brendle-log}, in our more general setting we need to introduce  the right function space that works well to ensure that the calculations involved in our proofs are formally correct. 
			In order to do that, we denote by  $W^{1,2}(\Sigma,d{\rm vol}_\Sigma)=\{f\in L^2(\Sigma,d{\rm vol}_\Sigma):\nabla^\Sigma f\in L^2(\Sigma,d{\rm vol}_\Sigma)\}$ the usual space of Sobolev functions on $\Sigma$. Furthermore, since we should also control  the term $|H|f$ for  $f\in W^{1,2}(\Sigma,d{\rm vol}_\Sigma)$, arising from the presence of the mean curvature, we introduce the function space
	$$W_H^{1,2}(\Sigma,d{\rm vol}_\Sigma)=\{f\in W^{1,2}(\Sigma,d{\rm vol}_\Sigma): |H|f\in L^2(\Sigma,d{\rm vol}_\Sigma)\}.$$  
%As we shall see in the sequel, it is precisely this framework that  involves the mean curvature vector $H$ is the one that will suit our purposes. 
As we will see in the sequel, this is precisely the appropriate framework,  that meets our goals.
When $H$ is bounded on $\Sigma$ (e.g., when $\Sigma$ is compact, or $\Sigma$ is a minimal submanifold), we clearly have $W_H^{1,2}(\Sigma,d{\rm vol}_\Sigma)=W^{1,2}(\Sigma,d{\rm vol}_\Sigma)$.

	Our first main result reads as follows: 
	
%	 After the influential works of Michael and Simon \cite{MS} and Allard \cite{Allard}, a huge interest is devoted to  Sobolev-type inequalities on Euclidean submanifolds.  

		\begin{theorem}\label{main-theorem-Brendle} Let $n\geq2$ and $m\geq 1$ be integers and  $\Sigma$ be a complete $n$-dimensional submanifold  of $\mathbb R^{n+m}$ without boundary.  Then for every $f\in W_H^{1,2}(\Sigma,d{\rm vol}_\Sigma)$ with $\int_\Sigma f^2d{\rm vol}_\Sigma=1$ we have 
		\begin{equation}\label{main-inequality-Brendle}
			\int_\Sigma f^2\log f^2 d{\rm vol}_\Sigma\leq \frac{n}{2}\log\left(\frac{2}{\pi e n }			\int_\Sigma \left(|\nabla^\Sigma f|^2+\frac{1}{4} |H|^2f^2\right)d{\rm vol}_\Sigma\right).
		\end{equation} 
		 Moreover, inequality  \eqref{main-inequality-Brendle} is sharp, and equality holds  for some function  $f\in W_H^{1,2}(\Sigma,d{\rm vol}_\Sigma)$ if and only if $\Sigma $ is isometric to the Euclidean space $\mathbb R^n$ and 
		 $f(x)=\left(\frac{\alpha}{\pi}\right)^{\frac{n}{4}} e^{-\frac{\alpha}{2} |x-x_0|^2},$ $x\in \mathbb R^n,$ where $\alpha>0$ and $x_0\in\mathbb R^n$.
%		
%		 $f(x)=(2\alpha)^{-\frac{n}{2}} e^{-\alpha |x|^2}$, $x\in \mathbb R^n$.
	\end{theorem}
	
%	Clearly, \eqref{main-inequality-Brendle} is meaningful whenever  $\int_\Sigma f^2|H|^2<+\infty;$ this fact surely occurs for every $f\in C_0^\infty(\Sigma)$, for instance, when $\Sigma\subset \mathbb R^{n+m}$ is compact.
Note that the right hand side of \eqref{main-inequality-Brendle} is finite as  $f\in W_H^{1,2}(\Sigma,d{\rm vol}_\Sigma)$. As already mentioned, the proof of Theorem \ref{main-theorem-Brendle} uses optimal transport theory. This idea goes back to the  paper of Cordero-Erausquin, Nazaret and Villani \cite{CENV}, where the authors used the method of optimal transport to prove sharp Sobolev inequalities in the Euclidean space. In order to follow this strategy we have to develop the optimal transport method for Euclidean submanifolds. Results in this direction have been recently developed by Wang \cite{Wang}, based on the papers by McCann \cite{McCann}, McCann and Pass \cite{McCann-Pass} and Cordero-Erausquin, McCann and Schmuckenschl\"ager \cite{CEMS}; however, we cannot apply Wang's results as they are valid only for compactly supported measures. Instead, we shall prove a stronger statement (see Theorem \ref{Monge-Ampere-main}), based on McCann \cite{McCann-Duke}, that works for not necessarily compactly supported measures as well. A similar argument can be also found in Figalli and Gigli \cite{Figalli-Gigli} on non-compact manifolds.

	 A simple consequence of  Theorem \ref{main-theorem-Brendle} is the parametric  $L^2$-logarithmic-Sobolev inequality: 
	
	\begin{corollary}\label{corollary-Brendle}
		Under the same assumptions as in Theorem \ref{main-theorem-Brendle}, 
		for any fixed $\alpha>0$ 
		and $f\in W_H^{1,2}(\Sigma,d{\rm vol}_\Sigma)$ with $\int_\Sigma f^2d{\rm vol}_\Sigma=1,$ one has
		 \begin{equation}\label{parametric}
		 	\int_\Sigma f^2\log f^2d{\rm vol}_\Sigma \leq -n+\frac{n}{2}\log\left(\frac{\alpha}{\pi}\right)+\frac{1}{\alpha}\int_\Sigma \left(|\nabla^\Sigma f|^2+\frac{1}{4} |H|^2f^2\right)d{\rm vol}_\Sigma.
		  \end{equation}
	  Moreover, equality holds in \eqref{parametric} for some  $f\in W_H^{1,2}(\Sigma,d{\rm vol}_\Sigma)$ if and only if $\Sigma $ is isometric to the Euclidean space $\mathbb R^n$ and  there exists  $x_0\in\mathbb R^n$ such that  $f(x)=(\frac{\alpha}{\pi})^{\frac{n}{4}} e^{-\frac{\alpha}{2} |x-x_0|^2}$, $x\in \mathbb R^n$.
	\end{corollary}
	 We notice that Brendle's sharp inequality  \eqref{Brendle-initial-inequality} directly follows from the first part of Corollary \ref{corollary-Brendle} 
by choosing $\alpha:=\frac{1}{4}$ and $f:=g^\frac{1}{2}/(\int_\Sigma g)^{1/2}$ for any smooth function $g>0$ on the compact submanifold $\Sigma\subset \mathbb R^{n+m}$ without boundary. More surprisingly, it turns out that the \textit{parametric} version \eqref{parametric} of the $L^2$-logarithmic-Sobolev inequality   is \textit{equivalent} to \eqref{main-inequality-Brendle}; see Remark \ref{remark-equivalence}.  
	
 Another consequence of   Corollary \ref{corollary-Brendle} can be stated for parametric Gaussian measures on $\Sigma$;  let $\alpha>0$ be fixed and 
 \begin{equation}\label{Gaussian-alpha-measure}
  d\gamma_\alpha(x)=\left(\frac\alpha\pi\right)^\frac{n}{2}e^{-\alpha|x|^2}d{\rm vol}_\Sigma(x),\ \ x\in \Sigma.
 \end{equation}
To state this, we consider the weighted Sobolev space defined by  
$$W_{H,1}^{1,2}(\Sigma,d\gamma_\alpha)=\{\varphi\in W_H^{1,2}(\Sigma,d\gamma_\alpha):|x|\varphi\in L^2(\Sigma,d\gamma_\alpha)\}.$$
 \begin{corollary}\label{corollary-gauss}
 	Under the same assumptions as in Theorem \ref{main-theorem-Brendle}, 
 	for any fixed $\alpha>0$ 
 	and $\varphi\in W_{H,1}^{1,2}(\Sigma,d\gamma_\alpha)$ with $\int_\Sigma \varphi^2d\gamma_\alpha=1,$ one has
 	\begin{equation}\label{parametric-gauss}
 		\int_\Sigma \varphi^2\log\varphi^2 d\gamma_\alpha\leq \frac{1}{\alpha}
 		\int_\Sigma |\nabla^\Sigma \varphi|^2d\gamma_\alpha+\frac{1}{4\alpha}\int_\Sigma \varphi^2|H+2\alpha x^\perp|^2d\gamma_\alpha.
 	\end{equation}
Moreover, equality holds in 
 \eqref{parametric-gauss} for some $\varphi\in W_{H,1}^{1,2}(\Sigma,d\gamma_\alpha)$ if and only if $\Sigma$ is isometric to the Euclidean space $\mathbb R^n$ and there exists  $x_0\in \mathbb R^n$ such that  $\varphi(x)=e^{\alpha \langle x,x_0\rangle -\frac{\alpha}{2}|x_0|^2}$, $x\in \mathbb R^n$.
 \end{corollary}
 Let us observe that the characterization of the equality case in the above statement answers a question raised in \cite{Wang-Xia-Zhang}.
Note that condition $\varphi\in W_{H,1}^{1,2}(\Sigma,d\gamma_\alpha)$ is sufficient to guaranty that the right hand side of \eqref{parametric-gauss} is finite. In fact, it seems that the space  $W_{H,1}^{1,2}(\Sigma,d\gamma_\alpha)$ is the optimal choice to prove Corollary \ref{corollary-gauss}, which is based on Corollary  \ref{corollary-Brendle} and a  divergence theorem on not necessarily compact domains (that is new even in the Euclidean setting), see Lemma \ref{divergencia};  a detailed discussion on the optimality of $W_{H,1}^{1,2}(\Sigma,d\gamma_\alpha)$ can be found in Remark \ref{remark-space}.

 As before, if $\alpha:=\frac{1}{4}$ and $\varphi:=g^\frac{1}{2}/(\int_\Sigma gd\gamma_\alpha)^{1/2}$ for any smooth function $g>0$ on the compact submanifold $\Sigma\subset \mathbb R^{n+m}$ without boundary, the first part of Corollary \ref{corollary-gauss} provides precisely the same conclusion as  Brendle \cite[Corollary 2]{Brendle-log}. 

%Considering the case of equality in Theorem \ref{main-theorem-Brendle} we find the following rigidity type statement that is similar formally similar to the rigidity result of the authors \cite{BK} for the equality case Borell-Brascamp-Lieb inequality on curved spaces. 
%
%\begin{theorem} \label{theorem-equality-case} 
%Let $n\geq2$ and $m\geq 1$ be integers and  $\Sigma$ be a complete $n$-dimensional submanifold of $\mathbb R^{n+m}$.  Assume, that there exists a non-zero function $f\in W^{1,2}(\Sigma)$  with $\int_\Sigma f^2=1$ satisfying  
%		\begin{equation}\label{main-equality-Brendle}
%			\int_\Sigma f^2\log f^2 d{\rm vol}_\Sigma =  \frac{n}{2}\log\left(\frac{2}{\pi e n }			\int_\Sigma \left(|\nabla^\Sigma f|^2+\frac{1}{4} f^2|H|^2\right)d{\rm vol}_\Sigma\right).
%		\end{equation} 
%	  Then $\Sigma $ is isometric to the Euclidean space $\mathbb R^n$ and $f(x)=(2\alpha)^{-\frac{n}{2}} e^{-\alpha |x-x_0|^2}$, $x\in \mathbb R^n$.
%\end{theorem}

It is worthwhile to point out that Theorem \ref{main-theorem-Brendle} (as well as Corollaries \ref{corollary-Brendle} and \ref{corollary-gauss}) is a particular case of a more general result that holds for any $p\geq 2$, see Theorem \ref{theorem-general} in \S\ref{Section-4}.  The proof of  the latter result is much more technical and the constants involved in its statement are codimension-dependent unless $p=2$. However, we can formulate here another statement for $p\geq 2$ with an elegant formula of the constant that is valid for \textit{minimal submanifolds}. In the proof we have to combine the  optimal transport method with a refined monotonicity property of the Gamma function. 
%In this way, we obtain a codimension-free $L^p$-logarithmic-Sobolev inequality for every $p\geq 2$ on {\it minimal submanifolds} of any codimension.  
To state this result we denote for $p>1$ the Sobolev space $W^{1,p}(\Sigma,d{\rm vol}_\Sigma)=\{f\in L^p(\Sigma,d{\rm vol}_\Sigma):\nabla^\Sigma f\in L^p(\Sigma,d{\rm vol}_\Sigma)\}$. Using this notation we can state:  
	
	\begin{theorem}\label{main-theorem} Let $n\geq2$ and $m\geq 1$ be integers and $p\geq 2$.
		Let $\Sigma$ be a complete $n$-dimensional  minimal submanifold of $\mathbb R^{n+m}$.  Then for every $f\in W^{1,p}(\Sigma,d{\rm vol}_\Sigma)$ with $\int_\Sigma |f|^pd{\rm vol}_\Sigma=1$ we have 
		\begin{equation}\label{main-inequality}
			\int_\Sigma |f|^p\log |f|^p d{\rm vol}_\Sigma\leq \frac{n}{p}\log\left( \left(\frac{p^2}{2\pi e n}\right)^{\frac{p}{2}}\int_\Sigma |\nabla^\Sigma f|^pd{\rm vol}_\Sigma\right).
		\end{equation} 
	Moreover, \eqref{main-inequality} is sharp for $p=2$ 	 and the equality characterization reads as in Theorem \ref{main-theorem-Brendle}. 
	\end{theorem}

Compared to the result of Del Pino and Dolbeault \cite{delPinoDolbeault-2}, see \eqref{e-sharp-log-Sobolev} -- where $p<n$ is crucial -- in Theorem \ref{main-theorem} we do not need any restriction on $p$ with respect to the dimension $n$, except for the condition $p\geq 2$; a similar result can be found in Agueh,  Ghoussoub and Kang \cite[Corollary 3.3]{AGK} in the classical Euclidean setting. We note that for $p\in (1,2)$, our method does not provide a  codimension-free $L^p$-logarithmic-Sobolev inequality, even on minimal submanifolds. 

Beside the classical Euclidean spaces, our logarithmic-Sobolev inequalities can be applied on various non-compact submanifolds, as  catenoid, helicoid,  Enneper surface, 
Henneberg surface, Neouvius surface,
	Bour's  surface, all of them being minimal surfaces in $\mathbb R^3$. Higher codimensional examples of minimal surfaces include graphs of holomorphic functions. We refer to the  book of Osserman \cite{Osserman}  for more details on the theory of minimal surfaces. Further exotic, non-minimal structures are provided by $n$-dimensional Einstein submanifolds $\Sigma$ of $\mathbb R^{n+2}$, they being  warped products of the form $\Sigma=M_1\times_\varphi M_2$, where $M_1$ is an Einstein manifold (i.e., its Ricci  tensor is
	proportional to the Riemannian metric of $M_1$), etc. For further examples, see Dajczer and Tojeiro \cite{Dajczer-T}.

It is well known that logarithmic-Sobolev inequalities are strongly related  to hypercontractivity estimates, see e.g.\ \cite{BKT}, \cite{BobkovGL}, \cite{Fujita}, \cite{Gentil}. In the sequel, we present two applications concerning the hypercontractivity estimate for Hopf--Lax semigroups on submanifolds. For simplicity of exposition, we only consider the  case $p=2$.

Let $\Sigma$ be a complete $n$-dimensional submanifold of $\mathbb R^{n+m}$. For  $t>0$ and  $u:\Sigma\to \mathbb R$  we consider the Hopf--Lax semigroup given by the formula
%	\begin{equation}\label{inf-convolution-0}
	$$		{\bf Q}_{t}u(x):=\inf_{y\in \Sigma}\left\{u(y)+\frac{d^2_\Sigma(x,y)}{2t}\right\}, \ x\in \Sigma,
	$$
	where $d_\Sigma$ stands for the distance function on $\Sigma$ induced by the Riemannian metric coming from $\mathbb{R}^{n+m}$. 
	%	Clearly, the admissible $t>0$ and $u:X\to \mathbb R$ are fixed in such a way to ensure that 	${\bf Q}_{t}u(x)\in \mathbb R$;
	By convention, ${\bf Q}_{0}u=u.$ Since $\Sigma$ is not necessarily compact, we  have to impose certain assumptions on the function
	$u$ that is going to be precisely formulated in Section   \ref{section-5} (see also Balogh, Krist\'aly and Tripaldi \cite{BKT}). 
Then, if $0<a<b$, $H$ is bounded on $\Sigma$ and  $e^\frac{au}{2}\in W^{1,2}(\Sigma,d{\rm vol}_\Sigma)$, we state the following hypercontractivity estimate: 
	\begin{equation}\label{hypercontractivity-estimate-0}
\|e^{{\bf Q}_tu}\|_{L^b(\Sigma,d{\rm vol}_\Sigma)}\leq \|e^{u}\|_{L^a(\Sigma,d{\rm vol}_\Sigma)}\left(\frac{b-a}{2\pi t}\right)^{\frac{n}{2}\frac{b-a}{ab}}\left(\frac{a}{b}\right)^{\frac{n}{2}\frac{a+b}{ab}}e^{t\frac{\|H\|^2_{L^\infty(\Sigma)}}{6}\frac{a^2+ab+b^2}{a^2b^2}}.
	\end{equation}
Moreover, \eqref{hypercontractivity-estimate-0} is sharp in the sense that the factor $2\pi$ cannot be replaced by a greater constant, and we also characterize the equality case.  For the precise formulation of this result, including the additional assumptions on the function $u$, we refer to Theorem \ref{theorem-hypercontractivity}.

The second application is another hypercontractivity estimate for Hopf--Lax semigroups on \textit{self-similar shrinkers}, i.e., submanifolds $\Sigma\subset \mathbb R^{n+m}$ satisfying the equation  $H+\frac{x^\perp}{2}=0$. These objects appear as critical points of the Gaussian area $G(\Sigma)=(4\pi)^{-\frac{n}{2}}\int_\Sigma e^{-\frac{|x|^2}{4}}d{\rm vol}_\Sigma$, which are of crucial importance in the study of the mean curvature flow, see Huisken \cite{Huisken}. If $\Sigma\subset \mathbb R^{n+m}$ is a self-similar shrinker, and $d\gamma=d\gamma_{1/4}$ is the Gaussian measure (i.e., the measure from \eqref{Gaussian-alpha-measure} for $\alpha=\frac{1}{4}$), then for every $a,t>0$ and any smooth function $u$ in $\Sigma$, satisfying the sub-quadratic growth condition
\begin{equation}\label{growth}
	|u(x)| \leq C_1 + C_2|x|^{\theta} , \ x \in \Sigma,
\end{equation}
for some $C_1, C_2 >0 $ and $0 <\theta < 2 $, we prove the  hypercontractivity estimate:
	\begin{equation}\label{hypercontractivity-estimate-Gauss}
	\|e^{{\bf Q}_tu}\|_{L^{a+\frac{t}{2}}(\Sigma,d\gamma)}\leq \|e^{u}\|_{L^a(\Sigma,d\gamma)};
\end{equation}
moreover, the factor $2$ in  $L^{a+\frac{t}{2}}(\Sigma,d\gamma)$-norm cannot be replaced by a smaller constant. The equality in \eqref{hypercontractivity-estimate-Gauss} is also characterized, which happens if and only if $\Sigma$ is isometric to $\mathbb R^n$ and $u$ is an affine function; this characterization seems to be new also in the standard Euclidean setting. 
The precise statement is formulated as Theorem \ref{them-hyper-Gauss}, whose proof is based on Corollary \ref{corollary-gauss} and on the Euclidean volume growth property of self-similar shrinkers. The optimality of the range for $\theta$ in the growth condition \eqref{growth} is also discussed; see Remark \ref{remark-theta}.

%However, in \S\ref{section-3} we state a codimension-depending $L^p$-logarithmic-Sobolev inequality for a general %submanifold in the case when $p\geq 2$, see Theorem \ref{theorem-general}, which reduces to Theorem \ref{main-%theorem-Brendle} for $p=2$ and to Theorem \ref{main-theorem}  for minimal submanifolds. 

%\begin{theorem}
%	Equality holds if and only if  $\Sigma$ is isometric to $\mathbb R^n$ and up to isometry and translation, $f$ is any Gaussian of the form $f_\lambda(x)=(\lambda \pi^{-1})^{\frac{n}{2}} e^{-\lambda |x|^2}$, $x\in \mathbb R^n$, for every $\lambda>0$.
%\end{theorem}
	
%	We apply our logarithmic-Sobolev inequalities to derive Poincar\'e inequalities and hypercontractivity estimates for Hopf--Lax semigroups on submanifolds.
The paper is structured as follows. In Section \ref{section-2} we provide the statement and the proof of Theorem \ref{OMT-theorem-submanifold}, that will serve as the necessary background on the optimal transport theory in order to prove our main results. This section contains also a technical integration by parts inequality that is used in the sequel. In Section \ref{section-3} we prove our main results in the case $p=2$; namely, Theorem \ref{main-theorem-Brendle} and Corollaries \ref{corollary-Brendle} \&  \ref{corollary-gauss}. In Section \ref{Section-4} we consider the case $p\geq 2$, by  proving Theorem \ref{main-theorem} and Theorem \ref{theorem-general} (see below).  In Section \ref{section-5} we consider the hypercontractivity estimates, by proving \eqref{hypercontractivity-estimate-0} and \eqref{hypercontractivity-estimate-Gauss} together with discussion of the equality cases (see Theorems \ref{theorem-hypercontractivity} \& \ref{them-hyper-Gauss}). We conclude the paper with an Appendix containing the statements and proofs of some technical results about the Gamma function that we use in Section \ref{Section-4}.

	\section{Optimal transport on submanifolds}\label{section-2}
In this section we fix the notation  used in the paper and we state a Brenier-type  statement on optimal transport that we need to prove our main results.
 This is a similar statement to the one recently obtained by Wang \cite{Wang}, based on the papers of McCann \cite{McCann}, McCann and Pass \cite{McCann-Pass} and Cordero-Erausquin, McCann and Schmuckenschl\"ager \cite{CEMS}. The crucial difference between Wang's results and ours is that  we do not assume that the involved measures are compactly supported;  
 this will be important in the proof of Theorem \ref{main-theorem-Brendle}, in the characterization of the  equality case.  For a slightly different approach of the  application of optimal transport in proving Michael--Simon-type isoperimetric inequalities, we refer to Brendle and Eichmair \cite{Brendle-Eichmair}.
	
Let $\Sigma$ be a complete $n$-dimesional
 submanifold of $\mathbb R^{n+m}$, possibly with boundary. The gradient,  divergence, Laplacian  and Hessian associated to $\Sigma$ are denoted by $\nabla^\Sigma$, ${\rm div}^\Sigma$, $\Delta^\Sigma$ and $D_\Sigma^2$, respectively.  We use the notations of  $T\Sigma$ and $T^{\perp} \Sigma$ for the tangent and  normal bundles of $\Sigma$. For every $x\in \Sigma$ and $v\in \mathbb R^{n+m}$, we write $v=v^T+v^\perp$ as the orthogonal splitting of $v$ with $v^T\in T_x\Sigma$ and $v^\perp\in T_x^\perp\Sigma$, where $T_x\Sigma$ and $T_x^\perp\Sigma$ denote the tangent and normal spaces to $\Sigma$ at $x\in \Sigma,$ respectively. In addition, let ${\rm Pr}_{T^\perp\Sigma}:T^\perp{\Sigma}\to \Sigma$ be the bundle projection map.
 
 We recall that the second fundamental form of $\Sigma$ is a symmetric bilinear form 
 $II: T\Sigma\times T\Sigma \to T^{\perp}\Sigma$ such that  if $X,Y$ are two tangent vector fields and $V$ is a normal vector field to $\Sigma$, then $$\langle II(X,Y),V\rangle=\langle \overline D_XY,V\rangle=-\langle \overline D_XV,Y\rangle.$$
 Here, the notation $\overline D$ stands for the linear connection on the ambient Euclidean space $\mathbb R^{n+m}.$ Furthermore, the mean curvature vector $H$ is nothing but the trace of the second fundamental form $II$. The submanifold $\Sigma$ is \textit{minimal} if $H=0$ on every point of $\Sigma.$

 The cost function in $\mathbb R^{n+m}$ is given by 
 $c(x,y)=\frac{1}{2}|x-y|^2$. It turns out that we can replace the cost function $c$ by the scalar product cost $c(x,y) = -\langle x, y\rangle $ for $x,y \in \mathbb R^{n+m}$.  We also recall that if  $S\subseteq \mathbb R^{n+m}$ is a nonempty set, then the subgradient of a function $\psi:S\to \mathbb R\cup \{+\infty\}$ at the point $x\in D(\psi)$ in the effective domain $D(\psi)=\{x\in S:\psi(x)<+\infty\}$ is defined as  
 $$\partial \psi(x)=\{y\in \mathbb R^{n+m}:\psi(z)\geq \psi(x)+\langle y,z-x\rangle,\forall z\in S\}.$$

The following statement is the precise formulation of the optimal transport result that we shall use in the sequel:

\begin{theorem}\label{OMT-theorem-submanifold} Let $n\geq2$ and $m\geq 1$ be integers,  $\Sigma$ be a complete $n$-dimensional submanifold of $\mathbb R^{n+m}$, and $\Omega\subseteq \mathbb R^{n+m}$ be an open set. Let $\mu$ and $\nu$ be Borel probability measures on $\Sigma$ and $\Omega$ which are absolutely continuous with respect to $d{\rm vol}_\Sigma$ and $d\mathcal L^{n+m},$ respectively.  
	 Then there exist a measurable subset $A$ of the normal bundle $T^\perp \Sigma$ and a function $u:\Sigma \to \mathbb R \cup \{+\infty\}$ which is semiconvex in its effective domain and it is twice differentiable on the set ${\rm Pr}_{T^\perp \Sigma}(A)\subset \Sigma$, that will give rise to a map $\Phi:A\to \Omega$ given by 
	$$\Phi(x,v)=\nabla^\Sigma u(x)+v,\ \ (x,v)\in A,$$
	such that 
	\begin{itemize}
		\item[(i)]  Pythagorean's rule holds, i.e., $|\Phi(x,v)|^2=|\nabla^\Sigma u(x)|^2+|v|^2$ for every $(x,v)\in A;$
		\item[(ii)] $\mu({\rm Pr}_{T^\perp \Sigma}(A))=\nu(\Phi(A))=1;$ 
		\item[(iii)] the map $\Phi(x,\cdot):A\cap T_x^\perp \Sigma\to \partial u(x)\cap \Phi(A)$ is a bijection for every $x\in {\rm Pr}_{T^\perp \Sigma}(A);$
		\item[(iv)]   for every $(x,v)\in A$ the $n\times n$ matrix  $D_\Sigma^2u(x)-\langle II(x),v\rangle$ is symmetric and non-negative definite, and the determinant-trace inequality holds, i.e., 
		\begin{align}\label{determinant-trace-0}
			 {\rm det} D\Phi(x,v)&= {\rm det}[D_\Sigma^2u(x)-\langle II(x),v\rangle]\nonumber \\&\leq \left(\frac{\Delta_{\rm ac}^\Sigma u(x)-\langle H(x),v\rangle}{n}\right)^n,
		\end{align}
	 where $\Delta_{\rm ac}^\Sigma u$ stands for the absolute continuous part of the distributional Laplacian $\Delta_\mathcal D^\Sigma u;$
	\item[(v)] 	$\Delta_{\rm ac}^\Sigma u\leq  \Delta_{\mathcal D}^\Sigma u$ 
	in the sense of distributions, i.e.\ $$\ds\int_{\Sigma} f \Delta_{\rm ac}^\Sigma u \leq \ds\int_{\Sigma} f\Delta_{\mathcal D}^\Sigma u, $$  for any compactly supported smooth function  $f:\Sigma \to \mathbb R _+$.
		\item[(vi)] 
		if $F$ and $G$ are the density functions of the measures $\mu$ and $\nu$ with respect to the volume measures of $\Sigma$ and $\mathbb R^{n+m}$, then for  $\mu$-a.e.\ $x\in \Sigma$ we have the following integral version of  the 
		Monge--Amp\`ere equation 
		\begin{equation}\label{Monge-Ampere-main}
			F(x)=\int_{A\cap T_x^\perp \Sigma}G(\Phi(x,v)){\rm det} D\Phi(x,v)dv. 
		\end{equation}
		
	\end{itemize}
\end{theorem}

{\it Proof.} As we shall see below, the proof uses in an essential way the fact that $\Sigma$ is embedded into the Euclidean space $\mathbb R^{n+m}$. This fact allows to apply the results of the classical  convex analysis and obtain locally Lipschitz properties of the Kantorovich potential whose subgradient carries the support of a coupling of $\mu$ and $\nu$. 
However, in order to get good differentiability properties of the potential we need to restrict it to $\Sigma$ and apply the Alexandrov--Bangert--Rademacher regularity theorem to this restriction.

Let $\mu_\Sigma$ be  the measure on $\mathbb R^{n+m}$ defined by $$\mu_\Sigma(S)=\mu(S\cap \Sigma),$$
where $S\subseteq \mathbb R^{n+m}$ is any Borel set. It is clear that $\mu_\Sigma\in \mathcal P(\mathbb R^{n+m})$, i.e.,  $\mu_\Sigma$ is a Borel probability  measure on $\mathbb R^{n+m}$, having its support in  $\Sigma.$ Since $\mu_\Sigma,\nu\in \mathcal P(\mathbb R^{n+m})$, by McCann \cite[Theorem 6]{McCann-Duke}, there exists a  coupling $\gamma\in \mathcal P(\mathbb R^{n+m} \times \mathbb R^{n+m})$ such that $\mu_\Sigma$ and $\nu$  are its marginals and the support of $\gamma$, i.e., ${\rm supp}\gamma \subset \mathbb R^{n+m} \times \mathbb R^{n+m}$, is cyclically monotone with respect to the quadratic cost $c$ on $\mathbb R^{n+m}$. We notice here  that  McCann's result -- that we have just applied -- does not require the compactness of the supports of $\mu$ and $\nu$. We also note that this statement only guarantees the existence of a coupling with cyclically monotone support without its optimality. The next step is to observe that the cyclical monotonicity property of a subset of $\mathbb R^{n+m} \times \mathbb R^{n+m}$ with respect to the quadratic and scalar product costs are equivalent. In conclusion,  we find that the classical Rockafellar's theorem implies the existence of a convex function $\psi:\mathbb R^{n+m}\to \mathbb R\cup \{+\infty\}$ such that 
\begin{equation}\label{support_in_graph}
	{\rm supp}\,\gamma \subseteq {\rm Graph}\, \partial \psi=\{(x,y):x\in \mathbb R^{n+m}, y\in \partial \psi(x)\},
\end{equation}
where $\partial \psi(x)\subset \mathbb R^{n+m}$ is the usual subgradient of $\psi$ at the point $x\in \mathbb R^{n+m}$.

Since the measure $\mu_\Sigma$ is supported on $\Sigma$, we see that it does not vanish on  sets of Hausdorff dimension $n+m- 1$ and so we cannot claim that $\nabla \psi$ exists for $\mu$-a.e. point.  Instead, another approach is needed. The idea is to restrict the potential $\psi$ to $\Sigma$ and study its differentiability properties according to the method inspired by Cordero-Erausquin,  McCann and Schmuckenschl\"ager \cite{CEMS} and Wang \cite{Wang}.  

In order to carry out this plan, let us first consider  the Legendre transform $\psi^c:\mathbb R^{n+m}\to \mathbb R\cup \{+\infty\}$  of $\psi,$ defined by 
$$\psi^c(y)=\sup_{z\in \mathbb R^{n+m}}\{\langle z,y\rangle -\psi(z)\},\ y\in \mathbb R^{n+m},$$
which is also a convex function. By definition, it follows that 
\begin{equation}\label{1-equation-Legendre}
	\psi(z)+\psi^c(y) -\langle z,y\rangle \geq 0 \ \ \ {\rm for\ all}\ \ z,y\in \mathbb R^{n+m},
\end{equation}
and 
\begin{equation}\label{2-equation-Legendre}
(x,y)\in {\rm Graph}\, \partial \psi \Longleftrightarrow	\psi(x)+\psi^c(y) -\langle x,y\rangle = 0. 
\end{equation}
In particular, it follows that for every $(x,y)\in {\rm Graph}\, \partial \psi$ the values $\psi(x)$ and $\psi^c(y)$ are finite; moreover, due to  \eqref{support_in_graph}, the latter properties also hold for $\gamma$-a.e.\ $(x,y)\in \mathbb R^{n+m} \times \mathbb R^{n+m}$.  Since $\mu_\Sigma=\pi_1{_{\#}}\gamma$, it follows that $\psi(x)<+\infty$ for $\mu_\Sigma$-a.e.\ $x\in \mathbb R^{n+m}$;  thus, by the definition of $\mu_\Sigma$, we also have that $\psi(x)<+\infty$ for $\mu$-a.e.\ $x\in \Sigma$. 

Let us consider the effective domain of the convex potential $\psi$, i.e., 
$E:=D(\psi)=\{x\in \mathbb R^{n+m}:\psi(x)<+\infty\}$. We note that $\psi|_{E\cap \Sigma}$ is a locally Lipschitz function with respect to the Euclidean metric in $\mathbb R^{n+m}$. Since the Riemannian metric on $\Sigma$ and the classical Euclidean metric in $\mathbb R^{n+m}$ are locally equivalent,  $\psi|_{E\cap \Sigma}$ is a locally Lipschitz function with respect to the Riemannian metric on $\Sigma$.   Now, Rademacher's theorem implies that  $\psi|_{E\cap \Sigma}$ is differentiable $d{\rm vol}_\Sigma$-a.e.\ $x\in E\cap \Sigma$; clearly, differentiability in this context  is understood on the manifold $\Sigma$. 
%Let $E_k\subset E$, $k\in \mathbb N$, with $E_k\subset E_{k+1}$, be an exhaustion of $E=\cup_k E_k$ by compact and convex subsets. Note that $\psi|_{E_k}$ is locally Lipschitz with respect to the Euclidean metric in $\mathbb R^{n+m}$, $k\in \mathbb N$. Since the Riemannian metric on $\Sigma$ and the classical Euclidean metric in $\mathbb R^{n+m}$ are locally equivalent, it turns out that  $\psi|_{E_k\cap \Sigma}$ is locally Lipschitz with respect to Riemannian metric on $\Sigma$. In particular, $u:=\psi|_{\Sigma}$ is locally Lipschitz. 
In fact, from now on we disregard the fact that our function is possibly defined on a larger set; accordingly, we consider the function $u:\Sigma\to \mathbb R\cup \{+\infty\}$ be defined as $$u=\psi|_\Sigma.$$

In the sequel, we prove that $u$ is semiconvex on $E\cap \Sigma,$ i.e., for every $x_0\in E\cap \Sigma$ there exists a convex geodesic  ball $B_r^\Sigma(x_0)=\{y\in \Sigma:d_\Sigma(x_0,y)<r\}$ in $\Sigma$ for some $r>0$ and a smooth function $h:B_r^\Sigma(x_0)\to \mathbb R$ such that $u+h$ is geodesically convex on $B_r^\Sigma(x_0)\cap E$. According to \cite[Lemma 3.11]{CEMS}, it is enough to show that  the continuous function $u=\psi|_{E\cap \Sigma}$ verifies the property that  for every $x_0\in E\cap \Sigma$ there exists a neighborhood $S\subset E\cap \Sigma$ of $x_0$ and a constant $C>0$ such that for every $x\in S$ and $w\in T_x\Sigma$, 
\begin{equation}\label{semi-convexity-lemma}
	\liminf_{r\to 0}\frac{u(\exp_x^\Sigma(rw))+u(\exp_x^\Sigma(-rw))-2u(x)}{r^2}\geq -C,
\end{equation}
where $\exp_x^\Sigma:T_x\Sigma\to \Sigma$ is the usual exponential map on the Riemannian manifold $\Sigma$. In fact,  a closer inspection of the proof of \cite[Lemma 3.11]{CEMS} yields that it costs no generality to assume $|w|\leq 1$ in \eqref{semi-convexity-lemma}. Note that if $\exp_x^\Sigma(\pm rw)\notin E$ for small $r>0$, relation \eqref{semi-convexity-lemma} trivially holds.

Let us fix $x_0\in E\cap \Sigma$ and $r>0$ small enough (for instance, $B_{2r}^\Sigma(x_0)$ to be geodesically convex in $E\cap \Sigma$);  note that since $S:=\overline{B_{r}^\Sigma(x_0)}\subset B_{2r}^\Sigma(x_0)$ is a compact set in $E\cap \Sigma$,  its image $\partial u(S)$ is also compact in $\mathbb R^{n+m}$. In particular, since $u=\psi|_{\Sigma}$,  there exists $C_S>0$ such that  
\begin{equation}\label{basic-compacntness}
	|x-y|\leq C_S,\ \ \ \forall (x,y)\in {\rm supp}\, \gamma\subseteq {\rm Graph}\, \partial \psi \subseteq {\rm Graph}\, \partial u\ \ {\rm with}\ \ x\in S.
\end{equation}

On the one hand, for every $(x,y)\in {\rm supp}\, \gamma\subseteq {\rm Graph}\, \partial \psi$ with $x\in S$, and $w\in T_x\Sigma$ with $|w|\leq 1$,
by Wang \cite[Lemma A.1]{Wang} one can find $C>0$ such that 
\begin{align}\label{C-bound}
\nonumber	L:=&\liminf_{r\to 0}\frac{\langle y,\exp_x^\Sigma(rw)+\exp_x^\Sigma(-rw)-2x\rangle}{r^2}\\ \geq& -2|x-y||II(w)|-4|x-y|\coth(|x-y|)\nonumber \\&+ \liminf_{r\to 0}\frac{|\exp_x^\Sigma(rw)|^2+|\exp_x^\Sigma(-rw)|^2-2|x|^2}{r^2}\\
		\geq& -C.	\nonumber 
\end{align}
Indeed, the bounds for the first two terms in the right hand side simply follow by \eqref{basic-compacntness},  $|w|\leq 1$ and the continuity of the second fundamental form $II$, while the bound for the limit follows by basic properties of $\exp^\Sigma$, combined with the compactness of $S\subset \Sigma$.  

On the other hand, since ${\rm supp}\, \gamma\subseteq {\rm Graph}\, \partial \psi$,   $\mu_\Sigma=\pi_1{_{\#}}\gamma$ and $u=\psi|_{\Sigma}$, for every $(x,y)\in {\rm supp}\, \gamma$ with $x\in E\cap \Sigma$ we have 
$$u(\exp_x^\Sigma(rw))\geq u(x)+\langle y,\exp_x^\Sigma(rw)-x\rangle,$$ 
$$u(\exp_x^\Sigma(-rw))\geq u(x)+\langle y,\exp_x^\Sigma(-rw)-x\rangle,$$
for every $r>0$  and $w\in T_x\Sigma.$ It remains to add the above two inequalities and use the  limit estimate \eqref{C-bound} to conclude the proof of \eqref{semi-convexity-lemma}.

Due to Alexandrov--Bangert theorem (see \cite{McCann, Wang}), since $u$ is semiconvex in $E\cap \Sigma$, it turns out that it is also twice differentiable $d{\rm vol}_\Sigma$-a.e.\ on $E\cap \Sigma$. 
In particular, if $(x,y)\in {\rm supp}\,\gamma$ with $x\in E\cap \Sigma$, by \eqref{support_in_graph}, \eqref{1-equation-Legendre} and \eqref{2-equation-Legendre},   one has that the function $ z \to \langle z,y \rangle -u(z) , z\in \Sigma,$ has its maximum at $z=x$. Therefore, we obtain  for $d{\rm vol}_\Sigma$-a.e.\ $x \in E\cap \Sigma$ with 
$(x,y)\in {\rm supp}\, \gamma\subseteq {\rm Graph}\, \partial u$
 that
\begin{equation} \label{eq:sigma-grad-0}
	{\nabla^{\Sigma}}\big|_{z=x} \left(\langle z,y \rangle -u(z)\right) = 0,
\end{equation}
and 
\begin{equation} \label{eq:sigma-hess-0}
	D^2_{\Sigma}\big|{_{z=x}} \left(\langle z,y \rangle -u(z)\right) \leq 0.
\end{equation}
Note that relation \eqref{eq:sigma-grad-0} implies that $y^T=\nabla^{\Sigma}u(x)$;  therefore, there exists $v \in T_x^{\perp} \Sigma$ such that 
%$$ {\nabla^{\Sigma}}\big|_{z=x}(|z|^2/2 -\phi(z)) \in T_x \Sigma.$$
%Introducing the function $u: S \to \mathbb R$ defined by $u(x)= |x|^2/2 -\phi(x),$ we conclude that if $(x,y) $ is a point in the support of $\pi_{\mu, \nu}$ then there exists some $v \in T_x^{\perp} \Sigma$ such that 
\begin{equation} \label{eq:motivation-for-phi}
	y= \nabla^{\Sigma} u(x) +v \ \ \text{and}\ \  \ D^2_{\Sigma} u(x) - \langle II(x),v  \rangle \geq 0,
\end{equation} 
the second relation from \eqref{eq:motivation-for-phi} following by \eqref{eq:sigma-hess-0}.

The two combined conditions expressed in \eqref{eq:motivation-for-phi} serve as motivation to introduce the set 
$$A = \{(x,v) \in T^{\perp} \Sigma: \text{\eqref{eq:motivation-for-phi} holds and} \ (x,y)\in {\rm supp}\, \gamma \}, $$
and the mapping 
\begin{equation}\label{map-omt}
	 \Phi: A \to \mathbb R^{n+m}, \ \ \Phi(x,v)= \nabla^{\Sigma} u(x) +v .
\end{equation}
Then clearly  we have the inclusion
$ {\rm supp}\, \gamma \subseteq \{ (x, \Phi(x,v)): (x,v) \in A\}.$ Moreover, since 	 $\nabla^\Sigma u(x)\in T_x\Sigma$ and $v\in T_x^\perp\Sigma$ are orthogonal, one has the Pythagorean rule  
\begin{equation}\label{pitagorean}
	|\Phi(x,v)|^2=|\nabla^\Sigma u(x)|^2+|v|^2,
\end{equation}
which is property (i), while the first part of (iv) is precisely the second relation from \eqref{eq:motivation-for-phi}.  
%\textcolor{red}{By Let $II$ be the second fundamental form of $\Sigma$, i.e., if $X,Y$ are two tangent vector fields and $V$ is a normal vector to $\Sigma $, then $$\langle II(X,Y),V\rangle=\langle \overline D_XY,V\rangle=-\langle \overline D_XV,Y\rangle,$$
%where $\overline D$ stands for the linear connection on the ambient Euclidean space $\mathbb R^{n+m}.$ Moreover, the mean curvature vector $H$ is nothing but the trace of the second fundamental form $II$. The submanifold $\Sigma$ is minimal if $H=0$ on every point of $\Sigma.$}
In addition, in a properly chosen basis the differential of $\Phi(x,v)$ has  the  matrix representation 
$$D\Phi(x,v)=\left[\begin{matrix}
	D_\Sigma^2u(x)-\langle II(x),v\rangle & 0 \\
	* & {\rm id}_{m\times m} 
\end{matrix}\right],$$
where by \eqref{eq:motivation-for-phi}, it follows that the $n\times n$ matrix  $D_\Sigma^2u(x)-\langle II(x),v\rangle$ is symmetric and non-negative definite for every $(x,v)\in A$; see also Brendle \cite[Lemma 5]{Brendle-1} for a different but related reasoning. 
Therefore, by using the latter matrix relation and the arithmetic-geometric mean inequality, one has for every  $(x,v)\in A$  the determinant-trace inequality
$$
	0\leq {\rm det} D\Phi(x,v)= {\rm det}[D_\Sigma^2u(x)-\langle II(x),v\rangle]\leq \left(\frac{\Delta_{\rm ac}^\Sigma u(x)-\langle H(x),v\rangle}{n}\right)^n,
$$
which is precisely \eqref{determinant-trace-0}, concluding the proof of (iv). 
Similarly to \cite[Theorem 3.3]{Wang}, $A\subset T^{\perp} \Sigma$ is a measureable set, the map $\Phi:A \to \mathbb R^{n+m} $ is injective on $A$ and  
	\begin{equation}\label{A-set-definition}
			\nu(\Phi(A))=\mu({\rm Pr}_{T^\perp\Sigma}(A))=1,
	\end{equation}
which is property (ii). Moreover, (iii) follows in a similar manner as Corollary 3.4 in \cite{Wang}.  To prove (v), since $u$ is semiconvex on $E\cap \Sigma$, the singular part $\Delta_{\rm s}^\Sigma u$ of $\Delta_\mathcal D^\Sigma u$ is non-negative. Thus,  
$\Delta_{\rm ac}^\Sigma u\leq  \Delta_{\mathcal D}^\Sigma u$ 
in the sense of distributions, i.e.\ $\ds\int_{\Sigma} f \Delta_{\rm ac}^\Sigma u \leq \ds\int_{\Sigma} f\Delta_{\mathcal D}^\Sigma u $  for any compactly supported smooth function  $f:\Sigma \to \mathbb R _+$.

% one has
%$$
%\Delta_{\rm ac}^\Sigma u(x)\leq \Delta_{\rm ac}^\Sigma u(x)+\Delta_{\rm s}^\Sigma u(x)= \Delta_{\mathcal D}^\Sigma u(x).
%$$

%	 where ${\rm Pr}_{T^\perp\Sigma}:T^\perp{\Sigma}\to \Sigma$ is the bundle projection map.  
	
Let us note that up to this point we did not use the fact that $\nu$ is absolutely continuous w.r.t. $\mathcal{L}^{n+m}$, thus in particular claims (i)-(v) hold true in this generality. On the other hand, for  property (vi) this fact plays a crucial role. In fact, the proof of (vi) follows exactly in the same manner as  in  \cite[Theorem 3.5]{Wang} by using the disintegration of the coupling $\gamma$. The only difference is that since the support of $\nu$ is not necessarily compact we have to apply again Theorem 6 of \cite{McCann-Duke}. In this way we can conclude, as in the proof of Theorem 3.5 in \cite{Wang}, that if  $F$ and $G$ are the density functions of the measures $\mu$ and $\nu$ with respect to the volume measures of $\Sigma$ and $\mathbb R^{n+m}$, then for  $\mu$-a.e.\ $x\in \Sigma$ we have the following integral version of  the 
Monge--Amp\`ere equation 
\begin{equation}\label{M-ampere}	F(x)=\int_{A_x}G(\Phi(x,v)){\rm det} D\Phi(x,v)dv, 
\end{equation}
where $A_x=A\cap T_x^\perp\Sigma$ for $x\in {\rm Pr}_{T^\perp\Sigma}(A)$, which proves (vi). \hfill $\square$

\begin{remark}\rm A similar change of variable formula to Monge--Amp\`ere equation \eqref{M-ampere} can be found in the paper by McCann and Pass \cite{McCann-Pass}, where the authors developed the optimal mass transport theory from higher dimension to lower dimension with general cost function. 
\end{remark}

%
%Useful byproducts of the above proof can be summarized as follows: 
%
%\begin{corollary}\label{corollary-OMT} Let us consider the set $A\subset T^\perp \Sigma$ and the functions $u:\Sigma\to \mathbb R \cup \{+\infty\}$ and $\Phi:A\to \mathbb R$ from Theorem \ref{OMT-theorem-submanifold}. Then the following inequalities hold$:$ 
%	\begin{itemize}
%		\item[(i)] {\rm (Determinant-trace inequality)} \begin{equation}\label{determinant-trace}
%			{\rm det} D\Phi(x,v)\leq \left(\frac{\Delta_{\rm ac}^\Sigma u(x)-\langle H(x),v\rangle}{n}\right)^n, \ \forall (x,v)\in A,
%		\end{equation}
%	
%	\end{itemize} 
%\end{corollary}
%
%{\it Proof.} 
%(i) 
%
%
%
%(ii)  
%\hfill $\square$\\

%	\subsection{An integration by parts inequality}

We conclude this section by proving an \textit{integration by parts inequality}, that plays a crucial role in the proof of our main results. For a Euclidean version of this result we refer to Cordero-Erausquin,  Nazaret and Villani \cite[Lemma 7]{CENV}. 
In the sequel, let $p\geq 2$ and $q=\frac{p}{p-1}$ its conjugate. Moreover, let 	$W_H^{1,p}(\Sigma,d{\rm vol}_\Sigma)=\{f\in W^{1,p}(\Sigma,d{\rm vol}_\Sigma): |H|f\in L^p(\Sigma,d{\rm vol}_\Sigma)\}.$

\begin{proposition}\label{integration_by_parts} Let 
	$f\in W_H^{1,p}(\Sigma,d{\rm vol}_\Sigma)\cap C^\infty(\Sigma)$ and $G\in L^1(\mathbb R^{n+m})$ be non-negative functions with $$\int_\Sigma f^pd{\rm vol}_\Sigma=\int_{\mathbb R^{n+m}} G(y)dy=1\ \ and \ \ \int_{\mathbb R^{n+m}} |y|^qG(y)dy<+\infty.$$ Let  $u:\Sigma\to \mathbb R\cup \{+\infty\}$ be the  semiconvex function in its effective domain which is guaranteed in Theorem \ref{OMT-theorem-submanifold} for the measures $d\mu(x)=f^p(x)d{\rm vol}_\Sigma(x)$ and $d\nu(y)=G(y)dy$. Then the following inequality holds$:$
	\begin{equation}\label{integration-by-parts}
		\int_\Sigma f^p \Delta_{\rm ac}^\Sigma ud{\rm vol}_\Sigma \leq -p\int_\Sigma f^{p-1} \langle\nabla^\Sigma f, \nabla^\Sigma u\rangle d{\rm vol}_\Sigma.
	\end{equation}
\end{proposition}

{\it Proof.} On account of Theorem \ref{OMT-theorem-submanifold}, besides the function $u:\Sigma\to \mathbb R\cup \{+\infty\}$, we also consider the  measurable subset $A$ of the normal bundle $T^\perp \Sigma$ and the map $\Phi:A\to \mathbb R^{n+m}$ given by 
$\Phi(x,v)=\nabla^\Sigma u(x)+v,$ $(x,v)\in A.$

First, by Theorem \ref{OMT-theorem-submanifold}/(v) (after suitable approximation), one has that 
	\begin{equation}\label{int-parts-0}
		\int_\Sigma f^p \Delta_{\rm ac}^\Sigma u\, d{\rm vol}_\Sigma \leq \int_\Sigma f^p \Delta_{\mathcal D}^\Sigma u\, d{\rm vol}_\Sigma.
	\end{equation}

Let $x_0\in \Sigma$ and consider the sequence of functions $\{\chi_k\}_{k\in \mathbb N}\subset C_0^\infty(\Sigma)$ with $0\leq \chi_k\leq 1$, $\chi_k\equiv 1$ on $B_{k}^\Sigma(x_0)$, $\chi_k\equiv 0$ on $\Sigma\setminus B_{k+1}^\Sigma(x_0)$ and $|\nabla^\Sigma \chi_k|\leq M$ for some $M>0$ which is independent on $k\in \mathbb N$. Let $f_k=f\chi_k\in C_0^\infty(\Sigma)$. By Theorem \ref{OMT-theorem-submanifold}/(ii)  and the usual integration by parts we have that 
\begin{equation}\label{int-parts-1}
	\int_\Sigma f_k^p \Delta_{\mathcal D}^\Sigma u\, d{\rm vol}_\Sigma=  -p\int_\Sigma f_k^{p-1} \langle\nabla^\Sigma f_k, \nabla^\Sigma u\rangle d{\rm vol}_\Sigma.
\end{equation}
We intend to take the limit in \eqref{int-parts-1} when $k\to \infty$. First, we observe that 
\begin{equation}\label{111-}
	\lim_{k\to \infty}\int_\Sigma  f_k^{p-1}\langle \nabla^\Sigma u,\nabla^\Sigma f_k\rangle d{\rm vol}_\Sigma = \int_\Sigma f^{p-1} \langle\nabla^\Sigma u, \nabla^\Sigma f\rangle d{\rm vol}_\Sigma. 
\end{equation}
To prove \eqref{111-}, we have that $$|\nabla^\Sigma f_k|\leq |\nabla^\Sigma f|+Mf,$$ thus, by the fact that $0\leq f_k\leq f$, one has  $$ f_k^{p-1}|\langle \nabla^\Sigma u,\nabla^\Sigma f_k\rangle|\leq f^{p-1}|\nabla^\Sigma u| (|\nabla^\Sigma f|+Mf)$$ on $\Sigma.$ Note that the latter expression belongs to $L^{1}(\Sigma,d{\rm vol}_\Sigma).$  To see this, since $f\in  W_H^{1,p}(\Sigma,d{\rm vol}_\Sigma)$, one has first that $|\nabla^\Sigma f|+Mf \in L^{p}(\Sigma,d{\rm vol}_\Sigma).$ On the other hand, we also have 
$f^{p-1}|\nabla^\Sigma u|\in L^{q}(\Sigma,d{\rm vol}_\Sigma)$. Indeed, by the Pythagorean rule (see Theorem \ref{OMT-theorem-submanifold}/(i)), we have $|\nabla^\Sigma u(x)|\leq |\Phi(x,v)|$ for every $(x,v)\in A$. Therefore, by the Monge--Amp\`ere equation  (see Theorem \ref{OMT-theorem-submanifold}/(vi)) 	
\begin{equation}\label{MA-MA}
	f^p(x)=\int_{A_x}G(\Phi(x,v)){\rm det} D\Phi(x,v)dv, 
\end{equation}
and a change of variable, it follows that 
\begin{align*}
	\int_\Sigma f^{q(p-1)}|\nabla^\Sigma u|^q  d{\rm vol}_\Sigma=& \int_\Sigma f^{p}|\nabla^\Sigma u|^q  d{\rm vol}_\Sigma\\ 
\leq &\int_\Sigma \int_{A_x}|\Phi(x,v)|^q G(\Phi(x,v)){\rm det} D\Phi(x,v)dv d{\rm vol}_\Sigma(x)\\=&\int_{\mathbb R^{n+m}} |y|^qG(y)dy<+\infty.
\end{align*} 
Combining the above facts and since $f_k$ converges pointwisely to $f$ in $\Sigma$, the dominated convergence theorem implies \eqref{111-}. 

For the left hand side of \eqref{int-parts-1} we intend to apply Fatou's lemma, but we have no a priori information on the sign of $\Delta_{\mathcal D}^\Sigma u$. In order to overcome this problem,  for every $x\in {\rm Pr}_{T^\perp \Sigma}(A)$, let $v_x\in \overline A_x=\overline{A\cap T_x^\perp\Sigma}$ be such that $|v_x|=\inf_{v\in A_x}|v|.$ First, we observe that $x\mapsto f^p(x)\langle H(x),v_x \rangle$ belongs to $L^1(\Sigma,d{\rm vol}_\Sigma)$ and \begin{equation}\label{Lp-int-part}
	\lim_{k\to \infty}  \int_\Sigma f_k^p(x)\langle H(x),v_x \rangle d{\rm vol}_\Sigma(x) =\int_\Sigma f^p(x)\langle H(x),v_x \rangle d{\rm vol}_\Sigma(x).
\end{equation}
Indeed,  by the Pythagorean rule (see Theorem \ref{OMT-theorem-submanifold}/(i)), we have $|v_x|\leq |\Phi(x,v)|$ for every $(x,v)\in A$; thus,  a similar estimate as before shows that 
\begin{equation*}
	\int_\Sigma f^p(x)|v_x|^q d{\rm vol}_\Sigma(x)\leq \int_\Sigma f^{p}|\nabla^\Sigma u|^q  d{\rm vol}_\Sigma\leq\int_{\mathbb R^{n+m}} |y|^qG(y)dy<+\infty.
\end{equation*} 
Therefore, since  $f\in W^{1,p}_H(\Sigma,d{\rm vol}_\Sigma)$, H\"older's inequality implies that 
\begin{align*}
	I_1:=&\int_\Sigma f^p(x)|\langle H(x),v_x \rangle| d{\rm vol}_\Sigma(x)\\ \leq&  \int_\Sigma f^p(x)|H(x)||v_x| d{\rm vol}_\Sigma(x)\\ \leq& \left( \int_\Sigma |H|^pf^p d{\rm vol}_\Sigma\right)^\frac{1}{p}\left( \int_\Sigma f^p(x)|v_x|^q d{\rm vol}_\Sigma(x)\right)^\frac{1}{q}\\<&+\infty.
\end{align*}
Since $f_k$ converges pointwisely to $f$ in $\Sigma$, the dominated convergence theorem implies \eqref{Lp-int-part}.

%In a similar manner as above, by using now that $|\nabla^\Sigma u(x)|\leq |\Phi(x,v)|$ for every $(x,v)\in A$, 
%we also have that 
%$f_k^{p-1}\nabla^\Sigma u$ strongly converges to $f^{p-1}\nabla^\Sigma u$ in $L^{q}(\Sigma,d{\rm vol}_\Sigma).$ Moreover, we note that $\nabla^\Sigma f_k$ is bounded in $L^{p}(\Sigma,d{\rm vol}_\Sigma)$ (note that for every $k\in \mathbb N$ one has $|\nabla^\Sigma \chi_k|\leq M$ for some $M>0$) and $\nabla^\Sigma f_k$ converges  to $\nabla^\Sigma f$ in distributional sense; in particular, $\nabla^\Sigma f_k$  converges weakly to $\nabla^\Sigma f$ in $L^{p}(\Sigma,d{\rm vol}_\Sigma).$ 
%As a consequence, it yields that (up to a subsequence) 
%
%KIVETTEM

Note that $\Delta_\mathcal D^\Sigma u(x)-\langle H(x),v_x \rangle\geq\Delta_{\rm ac}^\Sigma u(x)-\langle H(x),v_x \rangle\geq 0$ for every $x\in {\rm Pr}_{T^\perp \Sigma}(A)$ (see Theorem \ref{OMT-theorem-submanifold}/(iv)\&(v)); thus, by Fatou's lemma one has that
 \begin{align}\label{Fatou-lemma}
	\nonumber I_2&:=	\int_\Sigma f^p(x) \left(\Delta_{\mathcal D}^\Sigma u(x)-\langle H(x),v_x \rangle\right) d{\rm vol}_\Sigma(x) \\&\leq \liminf_{k\to \infty}\int_\Sigma f_k^p(x) \left(\Delta_{\mathcal D}^\Sigma u(x)-\langle H(x),v_x \rangle\right) d{\rm vol}_\Sigma(x).
	\end{align}
%Since $u_{k,\varepsilon}^\delta$ converges to $u_{k,\varepsilon}$ in $W^{1,p}(\omega;E)$, one has that 
%$\nabla u_{k,\varepsilon}^\delta$ converges to $\nabla u_{k,\varepsilon}$ in $L^p(\omega;E)$  and $(u_{k,\varepsilon}^\delta)^{p-1}$ converges to $(u_{k,\varepsilon})^{p-1}$ in $L^{p'}(\omega;E)$ whenever $\delta\to 0$, the latter property coming from the property of superposition operators, see e.g.\ Willem \cite[Appendix A]{Willem}. In addition, since 
%$\nabla \phi$ is essentially bounded on ${\rm supp}(u_{k,\varepsilon}^\delta)$, we have that
%\[
%\int_E (u_{k,\varepsilon}^\delta)^{p-1}\nabla u_{k,\varepsilon}^\delta\cdot \nabla \phi \omega dx\rightarrow \int_E (u_{k,\varepsilon})^{p-1}\nabla u_{k,\varepsilon}\cdot \nabla \phi \omega dx \quad \text{as $\delta\to0$}.
%\] 
Now, taking the limit $k\to \infty$  in\vspace{0.2cm} \\
%\begin{eqnarray*}%\label{int-parts-2}
$ \displaystyle
	\int_\Sigma f_k^p(x) \left(\Delta_{\mathcal D}^\Sigma u(x)-\langle H(x),v_x \rangle\right) d{\rm vol}_\Sigma(x) $ $$ =  -p\int_\Sigma f_k^{p-1}(x) \langle\nabla^\Sigma f_k(x), \nabla^\Sigma u(x)\rangle d{\rm vol}_\Sigma(x)   -\int_\Sigma f_k^p(x)\langle H(x),v_x \rangle d{\rm vol}_\Sigma(x),$$
%\end{eqnarray*}
which in turn is equivalent to 
\eqref{int-parts-1}, 
and using relations \eqref{int-parts-0},  \eqref{111-}, \eqref{Lp-int-part} and \eqref{Fatou-lemma},  inequality \eqref{integration-by-parts} yields at once. 
\hfill $\square$

%An important ingredient in our argument is the sharp estimate of the Jacobian determinant ${\rm det} D\Phi(x,v)$ of the optimal transport map $\Phi$, appearing in \eqref{Monge-Ampere}. 

\section{Proof of  Theorem \ref{main-theorem-Brendle} and Corollaries \ref{corollary-Brendle}\&\ref{corollary-gauss}: the case $p= 2$}\label{section-3}
	
The present section is  devoted mainly to the proof of 
Theorem \ref{main-theorem-Brendle}, which  is divided into two main steps; we first deal with the proof of inequality  \eqref{main-inequality-Brendle}, and then we characterize the equality case. 
	
\subsection{Proof of inequality \eqref{main-inequality-Brendle}} 
	Let $f\in W_H^{1,2}(\Sigma,d{\rm vol}_\Sigma)$ be such that  $\int_\Sigma f^2d{\rm vol}_\Sigma=1$. In order to prove  
	\eqref{main-inequality-Brendle}, since $|\nabla^\Sigma |f|(x)|=|\nabla^\Sigma f(x)|$ for ${\rm vol}_\Sigma$-a.e.\ $x\in \Sigma$, we may assume that $f$ is non-negative;  moreover, by using density arguments, it is enough to consider $f\in W_H^{1,2}(\Sigma,d{\rm vol}_\Sigma)\cap C^\infty(\Sigma)$.

%	Let us assume that $f\in C_0^\infty(\Sigma)$ first and prove the inequality in the statement for such functions. The general case will then follow by the density of $C_0^\infty(\Sigma)$ in $W{1,2} (\Sigma)$. 
	
	Let $\alpha>0$ be arbitrarily fixed and 
%	take  $f\in C_0^\infty(\Sigma)$  with $\int_\Sigma f^2=1$. Let $S$ be the support of $f$, which is a compact subset with non-empty interior of $\Sigma$. 
	 for convenience, we introduce the constants $$M_0:=\int_{\mathbb R^{n+m}}e^{-\alpha|y|^2}dy\ \ {\rm and}\ \ M_1:=\int_{\mathbb R^{m}}e^{-\alpha |v|^2}dv.$$ 
%and for every $k\in \mathbb N\setminus \{0\}$, let 
%$$C_k:=\int_{\mathbb R^{n+m}}e^{-\alpha|y|^2}\mathbbm{1}_{B_k}(y)dy=\int_{B_k}e^{-\alpha|y|^2}dy,$$ where $ B_k$ is the closed ball in $\mathbb R^{n+m}$ with center at the origin and radius $k>0,$  and $\mathbbm{1}_E$ stands for the characteristic function of the set $E\subset \mathbb R^{n+m}.$ Note that all these numbers are finite and 
%\begin{equation}\label{CktoM0}
%\lim_{k\to \infty}C_k=M_0.
%\end{equation}
	In the sequel, we consider the  probability measures $$d\mu(x)=f^2(x) d{\rm vol}_\Sigma(x)\ \ {\rm  and}\ \  d\nu(y)=M_0^{-1}e^{-\alpha|y|^2} dy$$ on $\Sigma$ and $\mathbb R^{n+m}$, respectively. Then there exists a map $\Phi:A\to \mathbb R^{n+m}$ of the form 
	$$
	\Phi(x,v)=\nabla^\Sigma u(x)+v,
	$$
	where $u: \Sigma \to \mathbb R\cup \{+\infty\}$ and $A_x:=A\cap T_x^\perp\Sigma$ have the corresponding properties as in  Theorem \ref{OMT-theorem-submanifold}.
%	pushing-forward $(\Sigma ,\mu)$ into $(B_k,\nu_k)$  with respect to the cost function ${|\cdot|^2}/{2}$ on $\mathbb R^{n+m}$, having , i.e., 
%	
%	where $x\in \Sigma $ and  $ v\in  T^\perp_x\Sigma ,$
%	and $u_k:\mathbb R^{n+m}\to \mathbb R$ is a convex function.
	By the Monge--Amp\`ere equation \eqref{Monge-Ampere-main}, one has for ${\rm vol}_\Sigma$-a.e.\ $x\in \Sigma$ that
	 \begin{equation}\label{MA-first}
	 	f^2(x)=\frac{1}{M_0}\int_{A_x}e^{-\alpha |\Phi(x,v)|^2}{\rm  det}D\Phi(x,v)dv.
	 \end{equation}
	 The elementary inequality $y\leq e^{y-1}$ for  $y>0$ applied to $y = {\rm  det}^\frac{1}{n}D\Phi(x,v)$ and  the corresponding determinant-trace inequality (see \eqref{determinant-trace-0} in Theorem \ref{OMT-theorem-submanifold}/(iv))  imply for ${\rm vol}_\Sigma$-a.e.\ $x\in \Sigma$ and any $v\in A_x$ that  
	\begin{equation}\label{deter-trace}
		{\rm  det}D\Phi(x,v)\leq e^{n\left({\rm  det}^\frac{1}{n}D\Phi(x,v)-1\right)}\leq e^{\Delta_{\rm ac}^\Sigma u(x)-\langle H(x),v\rangle-n}.
	\end{equation}
	Therefore, by the  Pythagorean rule (see Theorem \ref{OMT-theorem-submanifold}/(i)), the Monge--Amp\`ere equation \eqref{Monge-Ampere-main} and by identifying  $T_x^{\perp } \Sigma $ with $\mathbb R^m$, we obtain for ${\rm vol}_\Sigma$-a.e.\ $x\in \Sigma$ that
	\begin{align}\label{A_x_inclusion}
		f^2(x) \leq& \frac{1}{M_0}e^{-\alpha |\nabla^\Sigma u(x)|^2+\Delta_{\rm ac}^\Sigma u(x)-n}\int_{\mathbb R^m}e^{-\alpha |v|^2-\langle H(x),v\rangle}dv\\  =&\nonumber\frac{1}{M_0}e^{-\alpha |\nabla^\Sigma u(x)|^2+\Delta_{\rm ac}^\Sigma u(x)-n+\frac{1}{4\alpha }|H(x)|^2}\int_{\mathbb R^m}e^{-\alpha \left|v+\frac{1}{2\alpha}H(x)\right|^2}dv\\ =&\frac{M_1}{M_0}e^{-\alpha |\nabla^\Sigma u(x)|^2+\Delta_{\rm ac}^\Sigma u(x)-n+\frac{1}{4\alpha }|H(x)|^2}.\nonumber
	\end{align}
In the last step of the above chain  we used a change of variables and the definition of the constant $M_1>0.$ Taking the logarithm of the latter inequality,   multiplying by $f^2$,   an integration over $\Sigma$  
	yields that
\begin{align}\label{equ-h-fss}
\nonumber 	\int_\Sigma f^2\log f^2  
	\leq& \log \frac{M_1}{M_0}-n-\alpha\int_{\Sigma} f^2|\nabla^\Sigma u|^2\\&+\int_{\Sigma} f^2 \Delta_{\rm ac}^\Sigma u+\frac{1}{4\alpha }\int_{\Sigma}|H|^2 f^2.
\end{align}
Hereafter, once we integrate on $\Sigma$ and  no confusion arises,  the  measure is considered to be $d{\rm vol}_\Sigma$.

 Note that $f|\nabla^\Sigma u|\in L^2(\Sigma,d{\rm vol}_\Sigma)$, thus the third term in the right hand side of \eqref{equ-h-fss} is well defined. Indeed, by the Pythagorean rule (see Theorem \ref{OMT-theorem-submanifold}/(i)), the Monge--Amp\`ere equation \eqref{Monge-Ampere-main} and a change of variable, it follows that 
\begin{align*}
	\int_\Sigma f^2|\nabla^\Sigma u|^2 d{\rm vol}_\Sigma \leq& \frac{1}{M_0}\int_\Sigma \int_{A_x}|\Phi(x,v)|^2 e^{-\alpha |\Phi(x,v)|^2}{\rm  det}D\Phi(x,v)dv d{\rm vol}_\Sigma(x)\\ \leq &\frac{1}{M_0}\int_{\mathbb R^{n+m}} |y|^2e^{-\alpha |y|^2}dy<+\infty.
\end{align*}

In what follows, we intend to integrate by parts in \eqref{equ-h-fss} the expression 
 involving the Laplacian operator.  This is the step that does not automatically apply since $f$ is not necessarily compactly supported and $\Delta_{\rm ac}^\Sigma$ is not the distributional Laplacian.  However, what we can do is to apply the {\it integration by parts inequality} formulated as Proposition \ref{integration_by_parts}; for $p=2$ relation \eqref{integration-by-parts} reads as
$$
\int_{\Sigma} f^2 \Delta_{\rm ac}^\Sigma u \leq -2\int_{\Sigma} f \langle \nabla^{\Sigma} f , \nabla^{\Sigma} u \rangle.$$
 Therefore, by \eqref{equ-h-fss} we obtain that 
	\begin{eqnarray*}
		\int_\Sigma f^2\log f^2 &\leq& \log \frac{M_1}{M_0}-n-\alpha\int_{\Sigma} f^2|\nabla^\Sigma u|^2-2\int_{\Sigma} f \langle \nabla^\Sigma f,\nabla^\Sigma u\rangle\\&&+\frac{1}{4\alpha }\int_\Sigma |H|^2f^2\\&=&\log \frac{M_1}{M_0}-n-\alpha\int_{\Sigma} \left|f\nabla^\Sigma u+\frac{1}{\alpha}\nabla^\Sigma f\right|^2+\frac{1}{\alpha}
		\int_{\Sigma} |\nabla^\Sigma f|^2\\&&+\frac{1}{4\alpha}\int_{\Sigma }|H|^2 f^2
		\\&\leq &\log \frac{M_1}{M_0}-n+\frac{1}{\alpha}
		\int_{\Sigma} |\nabla^\Sigma f|^2+\frac{1}{4\alpha}\int_{\Sigma}|H|^2 f^2.
	\end{eqnarray*}
Using formula \eqref{Gauss-integral}, it follows  that
	\begin{equation}\label{Brendle-version}
			\int_\Sigma f^2\log f^2 \leq -n+\frac{n}{2}\log\left(\frac{\alpha}{\pi}\right)+\frac{1}{\alpha}
		\left(	\int_\Sigma |\nabla^\Sigma f|^2+\frac{1}{4}\int_\Sigma |H|^2f^2\right).
	\end{equation}
	By minimizing the right hand side with respect to $\alpha>0$, we obtain that the minimum is realized for the value 
	\begin{equation}\label{alfa-choice}
			\alpha:=\frac{2}{n}\left(	\int_\Sigma |\nabla^\Sigma f|^2+\frac{1}{4}\int_\Sigma |H|^2f^2\right)>0.
	\end{equation}
	 Inserting this value of $\alpha$ into the right side of \eqref{Brendle-version}, we obtain that 
	$$\int_\Sigma f^2\log f^2 \leq \frac{n}{2}\log\left(\frac{2}{\pi e n }		\left(	\int_\Sigma |\nabla^\Sigma f|^2+\frac{1}{4}\int_\Sigma |H|^2f^2\right)\right),$$
	which is the desired inequality \eqref{main-inequality-Brendle}. 
	
	Due to the sharpness of \eqref{e-sharp-log-Sobolev}, and observing the equality of the two numbers 	$\frac{2}{\pi e n }=\mathcal L_{2,n},$ we obtain the sharpness of this constant in \eqref{main-inequality-Brendle}.
	
	\begin{remark}\rm 
		It is well-known that the  $L^p$-logarithmic-Sobolev inequality \eqref{e-sharp-log-Sobolev} in $\mathbb R^n$ 
can be obtained by  a limiting process from the sharp Gagliardo--Nirenberg--Sobolev inequality, see Del Pino and  Dolbeault \cite{delPinoDolbeault-JMPA} for $p=2$ and extended later for $p<n$ in \cite{delPinoDolbeault-2}.	A similar question arises in the setting of submanifolds for the inequality \eqref{main-inequality-Brendle}. As we noticed in the Introduction,  sharp Gagliardo--Nirenberg--Sobolev inequalities can be obtained on  minimal submanifolds  of
codimension at most two via rearrangement techniques (P\'olya--Szeg\H o inequality and co-area formula), which implies sharp $L^p$-logarithmic-Sobolev inequalities in the same geometric setting after a suitable limiting process as in \cite{delPinoDolbeault-JMPA, delPinoDolbeault-2}. However, on generic submanifolds this approach is not yet clear. 
In fact, in a recent work  we were able to prove only an \textit{asymptotically sharp}  $L^p$-Sobolev inequality on  minimal submanifolds, see Balogh, Krist\'aly and Mester \cite{BKM}. It is a challenging  problem to provide a sharp or even an asymptotically sharp
 Gagliardo--Nirenberg--Sobolev inequality on a general $n$-dimensional submanifold  of $\mathbb R^{n+m}$ for any codimension $m\in \mathbb N$.

%whether these type of
%		results can be obtained for Sobolev-Gagliardo-Nirenberg inequalities
%		
%	 addig a mi submanifold esetunkben nem igazan latjuk erre a lehetoseget, tekintettel arra, hogy a "klasszikus" Sobolev-egyenlotlenseget is csak aszimpotikus elesseggel tudtuk igazolni, es ide szepen becitalnank a J. London Math. Soc.-hoz bekuldott cikket, tehat az eles GNS-re a submanifold-ok eseten elegge neheznek tunik (ami majd elvileg implikalna egy eles log-Sobolev-et). Mit gondolsz? 
		
	\end{remark}
	
\subsection{Characterization of equality in \eqref{main-inequality-Brendle}}  In this subsection we are going to study the equality case in the statement of  Theorem \ref{main-theorem-Brendle}. The idea of the proof is as follows. Assuming the existence of an extremizer $f\in  W_H^{1,2}(\Sigma)$ in \eqref{main-inequality-Brendle}, we first  show that it is smooth, i.e., $f\in C^\infty(\Sigma)$, possibly with non-compact support. Then, since the proof of Theorem \ref{main-theorem-Brendle} is performed precisely for functions belonging to $W_H^{1,2}(\Sigma)\cap C^\infty(\Sigma),$  we can track back the inequalities and identify in this way the equalities. The proof is divided into two steps.

%Assuming the existence of an extremal function in \eqref{main-inequality-Brendle}, we show first that this function is smooth, possibly with non-compact support. Next we show that using our smooth extremal function in the proof of Theorem \ref{main-theorem-Brendle}, we get equalities in the chain of inequalities within the proof,  

		\subsubsection{\it Regularity of extremizer.} Let us assume that there exists an   extremal function $f\in  W_H^{1,2}(\Sigma)$, i.e., equality holds in \eqref{main-inequality-Brendle}; as before,  we can assume without loss of generality that $f\geq 0$. By standard arguments from calculus of variations we conclude that $f$  is a weak solution of the elliptic PDE 
		\begin{equation} \label{euler-lagrange} 
			-\Delta^{\Sigma}f + |H|^2 f -c_1 f -c_2 f \log f= 0 \ \text{ in} \ \Sigma,
		\end{equation} 
		for some $c_1, c_2 \in \mathbb R$. Here, we set by continuity $f(x) \log f(x)=0 $  for $x \in \Sigma$ such that $f(x)=0$.

	We first prove that $f>0$ on $\Sigma;$ more precisely, for every compact  $K\subset \Sigma$ there exists $C_K>0$ such that 
		\begin{equation}\label{compact-positiv}
			f(x)\geq C_K \ \ {\rm for} \ \ {\rm vol}_\Sigma{\rm -a.e.} \ x\in K.
		\end{equation}
		 To see this,  we shall apply -- after a suitable adaptation --  the main result from the recent paper of Sirakov and Souplet \cite[Theorem 1.1]{Sirakov-Souplet}, which is recalled here for completeness:  
		
		 Let $\Omega \subset \mathbb R^{n}$ be  a domain of  the Euclidean space,  and let 
		$f: \Omega \to \mathbb R$ be a non-negative supersolution of 
		\begin{equation}\label{supersolution} 
			\mathcal L [f] \leq h(f) ,
		\end{equation}
		where $h: [0, \infty) \to \mathbb R$ is a continuous function with $h(0)=0$ and 
		\begin{equation}\label{Vazquez-limit}
			\limsup_{s \to 0} \frac{h(s)}{s(\log s)^{2}} < \infty. 
		\end{equation}
		Then, either $u\equiv 0$  in $\Omega$, or 
		$\text{essinf} _{B} f >0$ for every compact ball $B \subset \Omega$.

%		if  $\text{essinf} _{B} f =0$ for some ball $B \subset \subset \Omega$, then $ f$ must vanish identically on $\Omega$. 
		
		On the left side of the expression \eqref{supersolution},   $\mathcal L$ is an elliptic second order divergence operator that acts on $W^{1,2}(\Omega)$, having the form
		\begin{equation} \label{L}
			\mathcal L [u] := \text{div} (A(x) Du + b_{1}(x) u) + b_{2} (x)\cdot Du + c(x) u,
		\end{equation}
		 where the matrix-valued function $A$ belongs to $ L^{\infty}(\Omega) $ and there exist $0<\lambda < \Lambda$ with the property that $\lambda I \leq A(x) \leq \Lambda I$ for $x \in \Omega$, while the lower-order coefficients verify 		
		$b_{1}, b_{2} \in L_{\text{loc}}^{q}(\Omega)$ for some $q>n$, and $c \in  L_{\text{loc}}^{p}(\Omega)$ for some $p>n/2$.
		
		We intend to apply this result to our equation \eqref{euler-lagrange}; 
%		written in the form
%		$$ -\Delta^{\Sigma}f + (|H|^2(x) -c_1)  f = c_2 f \log f,  \ \text{in}  \ \Sigma.$$ 
			since it is understood in the setting of an embedded  Riemannian manifold, we shall work in local coordinate charts in order to reduce \eqref{euler-lagrange} to an equation in the Euclidean space. In local coordinates the expression of the Laplace-Beltrami operator $\Delta^{\Sigma}$ takes the form 
		$$  \Delta^{\Sigma}f = \frac{1}{\sqrt{|g|}} \partial_{i}(\sqrt{|g|}g^{ij}\partial_{j}f),$$
		where the usual Einstein summation convention is applied. 	
		By abuse of notation we understand the function $f$ to be our initial extremal function composed with local coordinate charts; thus we think of $f$ as  being defined on an open domain of $\mathbb R^{n}$.  Here we used the classical notation for $g_{ij}$ to be the components of the metric tensor matrix,  $g^{ij}$ are the components of its inverse, and $|g| := |\det g_{ij}|$.

		In order to bring $\Delta^{\Sigma}$ into the desired form $\mathcal L$ as in \eqref{L}, we notice that 
		\begin{eqnarray*} 
			\partial_{i}(g^{ij}\partial_{j}f) &=& \partial_{i} \left( \frac{1}{\sqrt{|g|}} \sqrt{|g|}g^{ij}\partial_{j}f \right) \\& =& \partial_{i} ( {1}/{\sqrt{|g|}}) \sqrt{|g|}g^{ij}\partial_{j}f + 
			\frac{1}{\sqrt{|g|}} \partial_{i} ( \sqrt{|g|}g^{ij}\partial_{j}f ) \\&=& -\partial_{i}(\log (\sqrt{|g|})) g^{ij}\partial_{j}f + \Delta^{\Sigma}f,
		\end{eqnarray*}
		which implies that 
		$$ \Delta^{\Sigma}f= \partial_{i}(g^{ij}\partial_{j}f) + \partial_{i}(\log (\sqrt{|g|})) g^{ij}\partial_{j}f  = \text{div} ( A(x) Df) + b(x)\cdot Df ,$$
		where $A = (g^{ij})_{ij}$ and the coordinates of  the vector field $b\in \mathbb R^{n}$ are given by $b_{j}= \partial_{i}(\log (\sqrt{|g|})) g^{ij}$ for $j=1, \ldots , n$. Thus, in the expression of $\mathcal L$ we may choose 
		$$A = (g^{ij})_{ij}, \ b_{1} = 0, \  b_{2} = b, \ \text{and} \ c(x) = -(|H(x)|^{2} -c_{1}), $$
%		such that 
%		$$ \mathcal L [f] =   \Delta^{\Sigma}f - (|H|^{2}(x) -c_{1})f = \text{div} (A(x) Df) + b(x)\cdot Df + c(x) f$$
		and the problem reduces to study the equation 
		$$ \mathcal L [f] = c_{2} f \log f  \ \ \ \ \text{in} \  \ \Omega .$$
		Here, $\Omega \subset \mathbb R^{n}$ is the domain of a coordinate chart whose image covers an open subset of $\Sigma$ and it is contained in a compact subset of $\Sigma$. 
		
		We notice that the assumptions of Theorem 1.1 of \cite{Sirakov-Souplet} are satisfied for the above chosen functions $A, b_{1}, b_{2}$ and $c$, respectively. Furthermore, the function $h(s) := c_{2} s \log s $, $s\geq 0,$ satisfies the condition  \eqref{Vazquez-limit}. Accordingly,   \cite[Theorem 1.1]{Sirakov-Souplet} applies;  since strict positivity of the solution is a local property, we obtain that our extremal function $f$  verifies  \eqref{compact-positiv}. 
		
		The positivity of $f$ combined with Sobolev embeddings gives the usual bootstrapping for the regularity, that yields $f\in C^{\infty}(\Sigma)$, see e.g.\ Nicolaescu \cite[Section 10.3]{Nicolaescu}.

		\subsubsection{\it Equality in \eqref{main-inequality-Brendle}.} Due to the previous argument, any  extremal function $f\in  W_H^{1,2}(\Sigma)$ in \eqref{main-inequality-Brendle} is positive (or, negative) on $\Sigma$,   belonging to $C^{\infty}(\Sigma)$. In particular, since the proof of inequality \eqref{main-inequality-Brendle} is provided precisely for this class of functions, the extremality of $f$ implies that inequalities have to be replaced by equalities. 
		
		 On the one hand, we should have equality in \eqref{deter-trace}, thus also in \eqref{determinant-trace-0}; observing that in the inequality $y \leq e^{y-1}$, $y>0$, equality holds only for $y=1$, we obtain for ${\rm vol}_\Sigma$-a.e.\ $x\in \Sigma$ and every $v \in A_x$ that 
		\begin{equation}\label{equal-1}
			{\rm  det}D\Phi(x,v)=1\ \ \ {\rm and}\ \  \Delta_{\rm ac}^\Sigma u(x)-\langle H(x),v\rangle-n=0.
		\end{equation}
		On the other hand, equality must also hold  in the inequality \eqref{A_x_inclusion} which in turn implies that $A_x = T_x^{\perp} \Sigma$  up to a ${\rm vol}_\Sigma$-null set. The latter fact combined with the second relation from \eqref{equal-1}  yields  that $H(x) = 0$ for ${\rm vol}_\Sigma$-a.e.\ $x\in \Sigma$. Moreover, we recall from \eqref{eq:motivation-for-phi} that for ${\rm vol}_\Sigma$-a.e.\ $x\in \Sigma$ and  for all $v\in A_x$ we have 
		\begin{equation} \label{fundamental form}
			\ D^2_{\Sigma} u(x) - \langle II(x),v  \rangle \geq 0.
		\end{equation} 
		Since $A_x = T_x^{\perp} \Sigma$ for ${\rm vol}_\Sigma$-a.e.\ $x\in \Sigma$, \eqref{fundamental form} is valid for  ${\rm vol}_\Sigma$-a.e.\ $x\in \Sigma$ and a.e.\ $v\in T_x^{\perp} \Sigma$; in particular,  not only the mean curvature $H(x)$, but also the entire second fundamental form $II(x)$ must vanish for ${\rm vol}_\Sigma$-a.e.\ $x\in \Sigma$. 
		This implies that the second fundamental form $II$ vanishes identically on $\Sigma$,  which means that $\Sigma$ is contained in an $n$-dimensional plane; for simplicity, we may assert that $\Sigma$ is isometric to a subset of $\mathbb R^n$.

 Since we should have equality in the integration by parts inequality \eqref{integration-by-parts} as well, the same should be true also in \eqref{int-parts-0}. In particular,  since  $f$ is positive and smooth on $\Sigma$, we necessarily have that the singular part of $\Delta_\mathcal D^\Sigma$  vanishes on $\Sigma$. Therefore, one has that $\nabla^\Sigma u\in W^{1,1}(\Sigma)$. By the equality in \eqref{determinant-trace-0}, it turns out that for some $\lambda>0$ one has $D^2_{\Sigma} u(x)=\lambda I_n$ for ${\rm vol}_\Sigma$-a.e.\ $x\in \Sigma$, where $I_n$ is the $(n\times n)$ unit matrix; in particular, by the first relation of \eqref{equal-1}, we necessarily have that $\lambda=1.$ Since $\Delta_{\rm s}^\Sigma u\equiv 0$ on $\Sigma\subseteq \mathbb R^n$, due to Figalli, Maggi and Pratelli \cite[Lemma A.2]{FMP}, there exists $x_0\in \Sigma$ such that $\nabla^\Sigma u(x)=x-x_0$ for a.e.\ $x\in \Sigma$. Therefore, $\Phi(x,v)=x+v-x_0$ for a.e.\ $x\in \Sigma$ and for any $v\in \mathbb R^m$.  Now, by the Monge-Amp\`ere equation \eqref{MA-first}, it turns out that for a.e.\ $x\in \Sigma$, 
	\begin{align}\label{Gaussian-identification}
	\nonumber	f^2(x)=&\frac{1}{M_0}\int_{A_x}e^{-\alpha |\Phi(x,v)|^2}{\rm  det}D\Phi(x,v)dv\\=&\nonumber\frac{1}{M_0}e^{-\alpha |x-x_0|^2}\int_{\mathbb R^m}e^{-\alpha |v|^2}dv\\=&\left(\frac{\alpha}{\pi}\right)^{\frac{n}{2}} e^{-{\alpha}|x-x_0|^2}.
	\end{align}
By \eqref{Gauss-integral}, \eqref{Gaussian-identification} and the fact $\int_\Sigma f^2=1,$ we have that
\begin{eqnarray*}
	1&=&\left(\frac{\alpha}{\pi}\right)^{\frac{n}{2}}\int_{\mathbb R^n} e^{-{\alpha}|x-x_0|^2}\\&=&\left(\frac{\alpha}{\pi}\right)^{\frac{n}{2}}\int_{\mathbb R^n\setminus \Sigma} e^{-{\alpha}|x-x_0|^2} + \left(\frac{\alpha}{\pi}\right)^{\frac{n}{2}}\int_{\Sigma} e^{-{\alpha}|x-x_0|^2}\\&=&\left(\frac{\alpha}{\pi}\right)^{\frac{n}{2}}\int_{\mathbb R^n\setminus \Sigma} e^{-{\alpha}|x-x_0|^2} + \int_{\Sigma} f^2 \\&=&\left(\frac{\alpha}{\pi}\right)^{\frac{n}{2}}\int_{\mathbb R^n\setminus \Sigma} e^{-{\alpha}|x-x_0|^2} + 1.
\end{eqnarray*}
In particular, it follows that 
$$\int_{\mathbb R^n\setminus \Sigma} e^{-{\alpha}|x-x_0|^2}=0,$$ which forces $\Sigma$ to be the entire space $\mathbb R^n$ (eventually up to a negligible set). In particular, by \eqref{Gaussian-identification},   $f$ is a Gaussian on $\mathbb R^n$ of the form provided by \eqref{Gaussian-identification}.

%			Falling into the setting of Del Pino and Dolbeault \cite{delPinoDolbeault-2}, see \eqref{e-sharp-log-Sobolev}, we conclude that the only extremizers in the $L^2$-logarithmic Sobolev inequality are the Gaussian functions belonging to the family  $f_\alpha(x)=(\frac{\alpha}{\pi})^{\frac{n}{4}} e^{-\frac{\alpha}{2} |x|^2}$, $x\in \mathbb R^n$, for every $\alpha>0$,  up to translations; see also Carlen \cite{Carlen} and the main result in \cite{BDK}.
		
		We notice that if we insert the  Gaussian from \eqref{Gaussian-identification} into  \eqref{alfa-choice}, we obtain an identity that holds for every $\alpha$; in particular, $\alpha>0$ can be arbitrarily chosen in \eqref{Gaussian-identification}.  
			\hfill $\square$
		
\begin{remark}\rm 
	The equality case in the Euclidean $L^2$-logarithmic-Sobolev inequality is in a perfect concordance with Carlen \cite{Carlen} and the main result of \cite{BDK}. 
\end{remark}

\textit{Proof of Corollary \ref{corollary-Brendle}.} Given $\alpha>0$ fixed, the required inequality \eqref{parametric} is precisely \eqref{Brendle-version}
in the proof of Theorem \ref{main-theorem-Brendle}. 
% Conversely, by the elementary inequality $\log(ey)\leq y$ for every $y>0$, inequality \eqref{main-inequality-Brendle} implies  \eqref{Brendle-version} for every $\alpha>0$.  
The equality case follows similarly as in Theorem \ref{main-theorem-Brendle}, obtaining that  $\Sigma$ is isometric to $\mathbb R^n$ and the 
only extremal function $f$ in \eqref{parametric} is provided by 
\begin{equation}\label{f-extremal}
	f^2(x)=\left(\frac{\alpha}{\pi}\right)^{\frac{n}{2}} e^{-{\alpha}|x-x_0|^2},\ \ x\in \mathbb R^n,
\end{equation}
for some $x_0\in \Sigma$.
  \hfill $\square$
%}

\begin{remark}\label{remark-equivalence}\rm By the proof of Theorem  \ref{main-theorem-Brendle}, it is clear that \eqref{main-inequality-Brendle} can be deduced by the parametric $L^2$-logarithmic-Sobolev inequality \eqref{Brendle-version} via a minimization argument. Alternatively, inequality \eqref{main-inequality-Brendle} together with the inequality $\log(ey)\leq y$ for every $y>0$ imply   \eqref{Brendle-version}. 
\end{remark}

%\textcolor{red}{
We now focus on the proof of Corollary  \ref{corollary-gauss}. To do this, we consider the  weighted Sovolev space $$W_{H,1}^{1,2}(\Sigma,d{\rm vol}_\Sigma)=\{f\in W_H^{1,2}(\Sigma,d{\rm vol}_\Sigma):|x|f\in L^2(\Sigma,d{\rm vol}_\Sigma)\},$$
and prove the following auxiliary result: 
\begin{lemma}\label{divergencia}
	For every $f\in W_{H,1}^{1,2}(\Sigma,d{\rm vol}_\Sigma)$, one has 
	\begin{equation}\label{difergencia-vanishes}
		\int_{\Sigma}{\rm div}^\Sigma(f^2 x^T)d{\rm vol}_\Sigma =0.
	\end{equation}
\end{lemma}
%}

	{\it Proof.}
Let us observe first that ${\rm div}^\Sigma(f^2 x^T) \in  L^1(\Sigma,d{\rm vol}_\Sigma)$ whenever  $f\in W_{H,1}^{1,2}(\Sigma,d{\rm vol}_\Sigma)$; accordingly,  the quantity 
$\int_{\Sigma}{\rm div}^\Sigma(f^2 x^T)d{\rm vol}_\Sigma $ is well-defined and  
$$\int_{\Sigma}{\rm div}^\Sigma(f^2 x^T)d{\rm vol}_\Sigma:= \lim_{ R \to \infty} \int_{\Sigma\cap B_R}{\rm div}^\Sigma(f^2 x^T)d{\rm vol}_\Sigma, $$
where $B_R$ is the usual ball in $\mathbb R^{n+m}$ 
with radius $R$ and centered at the origin. 
Indeed, we have that
$$|{\rm div}^\Sigma(f^2 x^T)| \leq 2 |f| \cdot |x| \cdot | \nabla^\Sigma f| + f^2(n+|H(x)|\cdot |x|),$$ thus H\"older's inequality combined with $f\in W_{H,1}^{1,2}(\Sigma,d{\rm vol}_\Sigma)$ provides the expected property.

Without loss of generality, we may assume that $f\in W_{H,1}^{1,2}(\Sigma,d{\rm vol}_\Sigma)\cap C^\infty(\Sigma)$ and $f$ is non-negative. Let $\{\chi_k\}_{k\in \mathbb N}\subset C_0^\infty(\Sigma)$ be the sequence of functions  from the proof of Proposition \ref{integration_by_parts}, and let $f_k=f\chi_k\geq 0$. It is clear that $f_k\in C_0^\infty(\Sigma)$ and the usual divergence theorem implies for every $k\in \mathbb N$ that 
\begin{equation}\label{1-diver}
	\int_{\Sigma}{\rm div}^\Sigma(f_k^2 x^T)d{\rm vol}_\Sigma =0.
\end{equation}
On the other hand, one has 
%since $\nabla^\Sigma f_k(x)\in T_x\Sigma $, it follows that $\langle \nabla^\Sigma f_k(x),x^\perp\rangle=0$, thus
\begin{eqnarray*}
	{\rm div}^\Sigma(f_k^2 x^T)&=& \langle \nabla^\Sigma (f_k^2)(x),x^T\rangle + f_k^2(x){\rm div}^\Sigma(x^T)\\&=&2 f_k(x) \langle \nabla^\Sigma f_k(x),x^T\rangle+f_k^2(x)(n+\langle H(x),x^\perp\rangle).
\end{eqnarray*}
Therefore, \eqref{1-diver} is equivalent to 
\begin{equation}\label{2-diver}
2 \int_{\Sigma}f_k(x) \langle \nabla^\Sigma f_k(x),x^T\rangle d{\rm vol}_\Sigma +	\int_{\Sigma}f_k^2(x)(n+\langle H(x),x^\perp\rangle)d{\rm vol}_\Sigma =0.
\end{equation}
We are going to take the limit in \eqref{2-diver}. Note that since $0\leq f_k\leq f,$ we have that
\begin{equation*}
	\ \left\{
	\begin{array}{lll}
	f_k(x) |\langle \nabla^\Sigma f_k(x),x^T\rangle|\leq f(x) (|\nabla^\Sigma f(x)|+Mf(x))|x^T|, \\
	\\
	f_k^2|n+\langle H(x),x^\perp\rangle|\leq f^2(n+|H(x)||x^\perp|),
	\end{array}%
	\right. 
\end{equation*}
where we used that for some $M>0$ one has $|\nabla^\Sigma \chi_k|\leq M$  for every $k\in \mathbb N$. Since $|x^T|\leq |x|$ and $|x^\perp|\leq |x|$, by H\"older's inequality and the fact that $f\in W_{H,1}^{1,2}(\Sigma,d{\rm vol}_\Sigma)$ one has that $f(x) (|\nabla^\Sigma f(x)|+Mf(x))|x^T|\in L^1(\Sigma,d{\rm vol}_\Sigma)$ and $ f^2(n+|H(x)||x^\perp|)\in L^1(\Sigma,d{\rm vol}_\Sigma).$
 By the pointwise convergence of $f_k$   to $f$ on $\Sigma$ and the dominated convergence theorem, if we take the limit in \eqref{2-diver} as $k\to \infty$, the above facts imply that $$	2\int_{\Sigma}f(x) \langle \nabla^\Sigma f(x),x^T\rangle d{\rm vol}_\Sigma+\int_{\Sigma}f^2(x)(n+\langle H(x),x^\perp\rangle)d{\rm vol}_\Sigma =0,$$
which is equivalent to \eqref{difergencia-vanishes}. \hfill $\square$

\begin{remark}\label{remark-space}\rm 
	Lemma \ref{divergencia} also seems to be unknown in the Euclidean case. In addition,  the choice of $W_{H,1}^{1,2}(\Sigma,d{\rm vol}_\Sigma)$ is optimal, i.e.,  in general we cannot consider 'larger' weighted Sobolev spaces than $W_{H,1}^{1,2}(\Sigma,d{\rm vol}_\Sigma)$ to guarantee the validity of \eqref{difergencia-vanishes}. To show this, we  restrict our attention to the Euclidean setting (thus $\Sigma=\mathbb R^n$ and $H=0$) and for 
%	that if. Indeed, 
every $\alpha\in [0,1]$, we consider $$W_{\alpha}^{1,2}(\mathbb R^n)=\{f\in W^{1,2}(\mathbb R^n):|x|^\alpha f\in L^2(\mathbb R^n)\}.$$
%It is clear that $W_{H,0}^{1,2}(\Sigma,d{\rm vol}_\Sigma)=W_{H}^{1,2}(\Sigma,d{\rm vol}_\Sigma)$. 
We claim that for every $\alpha\in [0,1)$,  one can find  $f\in W_{\alpha}^{1,2}(\mathbb R^n)$  such that \eqref{difergencia-vanishes} fails.  
 
 Indeed, let $f:\mathbb R^n\to [0,\infty)$ be a non-negative, radially symmetric function with $f(x)=h(|x|)$ and $h:[0,\infty)\to [0,\infty)$ is the saw-function defined by $h(s)=l_k-\frac{l_k}{\delta_k}|R_k-s|$ for  $s\in [R_k-\delta_k,R_k+\delta_k]$ and $0$ otherwise, where   $R_k$ diverges monotonically to $+\infty$ with  $R_k+\delta_k<R_{k+1}-\delta_{k+1}$, and $0<\delta_k,l_k\ll 1,$ $k\in \mathbb N.$ On the one hand, it is clear that $$\ds\int_{\mathbb R^n} f^2 \sim \sum_k l_k^2 R_k^{n-1}\delta_k, \b\ \ \ds\int_{\mathbb R^n} |\nabla f|^2 \sim \sum_k l_k^2 R_k^{n-1}\delta_k^{-1}$$ and $$\ds\int_{\mathbb R^n} |x|^{2\alpha}f^2 \sim \sum_k l_k^2 R_k^{2\alpha +n-1}\delta_k.$$ Now, choosing $R_k=2^k$, $l_k=R_k^{-n/2}$ and $\delta_k=k^\beta R_k^{-1}$ for any $\beta>1$, we obtain that $f\in W_{\alpha}^{1,2}(\mathbb R^n)$ whenever $\alpha\in [0,1)$. On the other hand, since $f$ is radially symmetric, and if  \eqref{difergencia-vanishes} holds, it follows that 
 $$0=	\int_{\mathbb R^n}{\rm div}(f^2 x)dx= \lim _{R \to \infty}\int_{\mathbb R^n \cap B_R}{\rm div}(f^2 x)dx =  n\omega_n\lim_{R\to \infty}h^2(R)R^n.$$ 
 However, $$\lim_{k\to \infty}h^2(R_k)R_k^n=\lim_{k\to \infty} l_k^2R_k^n=1,$$  a contradiction. As expected in view of Lemma \ref{divergencia}, for $\alpha=1$, one has $f\notin W_{1}^{1,2}(\mathbb R^n),$ since  $\ds\int_{\mathbb R^n} |x|^{2}f^2 \sim \sum_k l_k^2 R_k^{n+1}\delta_k=\sum_k k^\beta =+\infty.$ 
\end{remark}

%\textcolor{red}{
\textit{Proof of Corollary \ref{corollary-gauss}}. 
%As before, it is enough to consider compactly supported smooth functions.  
Let  $\alpha>0$ be fixed 
and $\varphi\in  W_{H,1}^{1,2}(\Sigma,d\gamma_\alpha)$ with $\int_\Sigma \varphi^2d\gamma_\alpha=1$. We consider the function $f:\Sigma \to \mathbb R$ defined by 
\begin{equation}\label{f-varphi}
	f(x)=\left(\frac\alpha\pi\right)^\frac{n}{4}e^{-\frac{\alpha}{2}|x|^2}\varphi(x),\ \ x\in \Sigma.
\end{equation}
 By definition, it is clear that $$f^2 d{\rm vol}_\Sigma=\varphi^2d\gamma_\alpha.$$ A straightforward computation shows that $f\in W_{H,1}^{1,2}(\Sigma,d{\rm vol}_\Sigma)$. Furthermore,    $\int_\Sigma f^2 d{\rm vol}_\Sigma=1$. Applying Corollary \ref{corollary-Brendle} we obtain 
\begin{align}\label{parametric-applied}
\nonumber	\int_\Sigma f^2\log f^2d{\rm vol}_\Sigma \leq & -n+\frac{n}{2}\log\left(\frac{\alpha}{\pi}\right)\\&+\frac{1}{\alpha}\int_\Sigma \left(|\nabla^\Sigma f|^2+\frac{1}{4} |H|^2f^2\right)d{\rm vol}_\Sigma.
\end{align}
Since $$\frac{\nabla^\Sigma f}{f}=\frac{\nabla^\Sigma \varphi}{\varphi}-\alpha x^T\ \  {\rm and}\ \ {\rm div}^\Sigma(x^T)=n+\langle H(x),x^\perp\rangle,$$
%$$d\gamma_\alpha(x)=\left(\frac\alpha\pi\right)^\frac{n}{2}e^{-\alpha|x|^2}d{\rm vol}_\Sigma(x),\ \ x\in \Sigma.$$
for every $x\in \Sigma$ (except for a ${\rm vol}_\Sigma$-null set), a simple computation yields that 
\begin{eqnarray*}
	\log f^2 +n+\frac{n}{2}\log\left(\frac{\pi}{\alpha}\right)-\frac{1}{\alpha}\frac{|\nabla^\Sigma f|^2}{f^2}-\frac{|H|^2}{4\alpha}-\frac{{\rm div}^\Sigma(f^2 x^T)}{f^2}
	\\
	=\log \varphi^2 -\frac{1}{\alpha}\frac{|\nabla^\Sigma \varphi|^2}{\varphi^2}-\frac{|H+2\alpha x^\perp|^2}{4\alpha}.
\end{eqnarray*}
Integrating the left hand side with respect to $f^2 d{\rm vol}_\Sigma$, while the right hand side  with respect to $\varphi^2d\gamma_\alpha$, and using Lemma \ref{divergencia}, i.e.,  $$\int_{\Sigma}{\rm div}^\Sigma(f^2 x^T)d{\rm vol}_\Sigma =0,$$ relation \eqref{parametric-applied} implies the required inequality \eqref{parametric-gauss}. 
%}

If equality holds in \eqref{parametric-gauss} for some $\varphi$, the integral form of the latter relation implies that it should be equality also in \eqref{parametric}  for $f$, where the functions $f$ and $\varphi$ are connected by relation \eqref{f-varphi}. Since  $f$ is given by \eqref{f-extremal}, the claim follows. 
\hfill $\square$

	\section{Proof of Theorem  \ref{main-theorem} and related results: the case $p\geq 2$} \label{Section-4}

\textit{Proof of Theorem  \ref{main-theorem}}. The strategy is to  follow the main idea of the proof of Theorem \ref{main-theorem-Brendle}. The difference is that we need to deal with the additional difficulty due to the fact that $p\geq 2$. 
 As before, we just consider functions $f\in W^{1,p}(\Sigma,d{\rm vol}_\Sigma)$ such that 
\begin{equation}\label{f-unit}
	\int_\Sigma |f|^p=1.
\end{equation}
%and  $S={\rm supp} f\subseteq \Sigma$. 

Let $\alpha>0$ be fixed. 
We consider the  probability measures $$d\mu(x)=|f|^p(x) d{\rm vol}_\Sigma(x) \ \ {\rm  and} \ \ d\nu(y)=M_0^{-1}e^{-\alpha|y|^q}dy,$$  respectively, where $q=\frac{p}{p-1}$ and as before
$$M_0:=\int_{\mathbb R^{n+m}}e^{-\alpha{|y|^q}}dy.$$
We apply Theorem \ref{OMT-theorem-submanifold} to conclude that there exists  $\Phi$ of the form
  $$
  \Phi(x,v)=\nabla^\Sigma u(x)+v,\ \ (x,v)\in A,
  $$
 	where $u:\Sigma\to \mathbb R \cup \{+\infty\}$, and $A_x:=A\cap T_x^\perp\Sigma$ have the corresponding properties.
By the Monge--Amp\`ere equation (see Theorem \ref{OMT-theorem-submanifold}/(vi)), one has for ${\rm vol}_\Sigma$-a.e.\ $x\in \Sigma$ that 
\begin{equation}\label{MA-p}
	|f|^p(x)=\frac{1}{M_0}\int_{A_x}e^{-\alpha|\Phi(x,v)|^q}{\rm det} D\Phi(x,v)dv
	.
\end{equation}
%where $A:=A_k$ and $A_x:=A_k\cap T_x^\perp\Sigma_0$ are  the corresponding sets from  \eqref{A-set-definition}.

To deal with the difficulty caused by $p \geq 2,$ we introduce an additional parameter $t\in (0,1)$ whose value  will be determined later. Since  $p\geq 2$, it follows that $1<q=\frac{p}{p-1}\leq 2;$ thus, the concavity of $|\cdot|^\frac{q}{2}$ implies that 
\begin{eqnarray*}
	|\Phi(x,v)|^q&=&\left(|\nabla^\Sigma u(x)|^2+|v|^2\right)^\frac{q}{2}=\left(t\frac{|\nabla^\Sigma  u(x)|^2}{t}+(1-t)\frac{|v|^2}{1-t}\right)^\frac{q}{2}\\&\geq& t^{1-\frac{q}{2}}\left|\nabla^\Sigma u(x)\right|^q+(1-t)^{1-\frac{q}{2}}|v|^q.
\end{eqnarray*}

For notational convenience, let 
\begin{equation}\label{M0-constant}
 M_t:=\int_{\mathbb R^{m}}e^{-\alpha{(1-t)^{1-\frac{q}{2}}}|v|^q}dv.
\end{equation}
By the above concavity estimate, relation \eqref{MA-p} and the determinant-trace inequality (see Theorem \ref{OMT-theorem-submanifold}/(iv)) together with the fact that $\Sigma$ is a minimal submanifold (thus, $H=0$ on $\Sigma$),    one has for ${\rm vol}_\Sigma$-a.e.\ $x\in \Sigma$ that
\begin{align}\label{log-int}
\nonumber	|f|^p(x)=&\frac{1}{M_0}\int_{A_x}e^{-\alpha{|\Phi(x,v)|^q}}{\rm det} D\Phi(x,v)dv\\ \leq& \frac{M_t}{M_0}e^{
		-\alpha{t^{1-\frac{q}{2}}}
		\left|\nabla^\Sigma u(x)\right|^q+\Delta_{\rm ac}^\Sigma u(x)-n},
\end{align}
where we used the notation from \eqref{M0-constant} and
$${\rm  det}D\Phi(x,v)\leq e^{n\left({\rm  det}^\frac{1}{n}D\Phi(x,v)-1\right)}\leq e^{\Delta_{\rm ac}^\Sigma u(x)-\langle H(x),v\rangle-n}=e^{\Delta_{\rm ac}^\Sigma u(x)-n}.$$ 
We proceed as in the  previous proof:  we  take the logarithm of \eqref{log-int}, multiply by $|f|^p\geq 0$, integrate over $\Sigma$ and use \eqref{f-unit} to obtain that 
\begin{align}\label{ineq-1}
	\int_\Sigma |f|^p\log |f|^p  \leq& \log \frac{M_t}{M_0} -n-\alpha{t^{1-\frac{q}{2}}}\int_{\Sigma} |f|^p\left|\nabla^\Sigma u\right|^q + \int_{\Sigma} |f|^p\Delta_{\rm ac}^\Sigma u .
\end{align}
%By using   Jensen's inequality, relation \eqref{distributional-ac}, an integration by parts and H\"older's inequality, we obtain that
%\begin{eqnarray*}
%	\int_{\Sigma_0} |f|^p(x)\log \left(\frac{\Delta_{\rm ac}^\Sigma u_k(x)}{n}\right)&\leq& \log\left(\frac{1}{n}\int_{\Sigma_0} |f|^p(x) \Delta_{\rm ac}^\Sigma u_k(x) \right)\\
%	&\leq& \log\left(\frac{1}{n}\int_{\Sigma_0} |f|^p(x) \Delta_\mathcal D^\Sigma u_k(x) \right)
%	\\&=&\log\left(-\frac{p}{n}\int_{\Sigma_0} |f|^{p-1}(x) \langle \nabla^\Sigma u_k(x),\nabla^{\Sigma} |f|(x)\rangle \right)\\&\leq &\log\left(\frac{p}{n}\left(\int_{\Sigma } |f|^{p}|\nabla^\Sigma u_k|^qdx \right)^\frac{1}{q}
%	\left (\int_{\Sigma } |\nabla^\Sigma f|^p\right)^\frac{1}{p}\right).
%\end{eqnarray*}	
By using  the \textit{integration by parts inequality} from Proposition \ref{integration_by_parts}  together with H\"older's inequality, we obtain by \eqref{ineq-1} that
\begin{align*}	\int_\Sigma |f|^p\log |f|^p \leq& 
	\log \frac{M_t}{M_0} -n-\alpha{t^{1-\frac{q}{2}}}\int_{\Sigma} |f|^p\left|\nabla^\Sigma u\right|^q\\& -p \int_{\Sigma} |f|^{p-1}\langle\nabla^\Sigma u , \nabla^\Sigma |f|\rangle\\ \leq  &\log \frac{M_t}{M_0} -n-\alpha{t^{1-\frac{q}{2}}}\int_{\Sigma} |f|^p\left|\nabla^\Sigma u \right|^q\\& +p \left(\int_{\Sigma} |\nabla^\Sigma f|^p\right)^\frac{1}{p} \left(\int_{\Sigma}|f|^p|\nabla^\Sigma u|^q\right)^\frac{1}{q}. 
\end{align*}
Let us note that in the previous case when $p=q=2$ we could handle the integral terms on the right side in a rather straightforward manner by completing the squared; this idea does not work in our present situation when $q\neq 2$. 
The trick is to consider the function 
$$ s\mapsto -\alpha{t^{1-\frac{q}{2}}} s +p \left(\int_{\Sigma } |\nabla^\Sigma f|^p\right)^\frac{1}{p}s^\frac{1}{q}, \ \ \ s\in [0,\infty),$$  
that appears naturally on the right side. 
An elementary calculation yields that this function has its maximum  at 
$s_0=\alpha^{-p}(p-1)^pt^{p(\frac{q}{2}-1)} \int_{\Sigma } |\nabla^\Sigma f|^p ,$ and its maximal value is equal to 
 $$\alpha^{1-p}(p-1)^{p-1}t^{(p-1)(\frac{q}{2}-1)} \int_{\Sigma } |\nabla^\Sigma f|^p.$$
  Replacing this maximal value in right side of  the above inequality we obtain that 
\begin{equation}\label{later-we-need}	\int_\Sigma |f|^p\log |f|^p \leq \log  \frac{M_t}{M_0} -n + \alpha^{1-p}(p-1)^{p-1}t^{(p-1)(\frac{q}{2}-1)} \int_{\Sigma} |\nabla^\Sigma f|^p.
\end{equation}
By using \eqref{Gauss-integral} from the Appendix,  it follows that
\begin{eqnarray*}
	\int_\Sigma |f|^p\log |f|^p &\leq&  \log\left(\frac{\omega_m \Gamma(\frac{m}{q}+1)}{\omega_{m+n} \Gamma(\frac{m+n}{q}+1)} (1-t)^{\frac{m}{q}(1-\frac{q}{2})}\right)-n \\&&+\frac{n}{q}\log \alpha+ \alpha^{1-p}(p-1)^{p-1}t^{(p-1)(\frac{q}{2}-1)} \int_{\Sigma} |\nabla^\Sigma f|^p.
\end{eqnarray*}
The latter estimate is valid for every $\alpha$; minimizing the right hand side with respect to the value of $\alpha$, a simple computation and manipulation of the terms yield that 
$$\int_\Sigma |f|^p\log |f|^p \leq \frac{n}{p}\log\left(\frac{p}{n}\left(\frac{p-1}{e}\right)^{p-1}K_{m,n,q}^\frac{p}{n}(t)\int_{\Sigma} |\nabla^\Sigma f|^p\right),$$
where $$K_{m,n,q}(t)=\frac{\omega_m \Gamma(\frac{m}{q}+1)}{\omega_{m+n} \Gamma(\frac{m+n}{q}+1)}  \left((1-t)^mt^n\right)^{\frac{1}{2}-\frac{1}{q}},\ \ t\in (0,1).$$ 
The above estimate holds for any $t\in (0,1)$. Using that $q \leq 2$ we note that $\frac{1}{2}-\frac{1}{q}\leq 0$ and thus we can minimize the function
$t\mapsto \left((1-t)^mt^n\right)^{\frac{1}{2}-\frac{1}{q}}$. Doing so, we obtain the minimum point at $t_0=\frac{n}{n+m}$;  that gives
%
%
%
% Since $q\leq 2,$ we choose $t_0=\frac{n}{n+m}$, which gives the best choice for the above expression;
%
%
%
%IDAIG. 
%
% A simple computation, based on \eqref{Gauss-integral} gives that
%$$\left(\frac{M_t}{M_0}\right)^\frac{p}{n}t^{(\frac{q}{2}-1)\frac{p}{q}}=q^{-\frac{p}{q}}\left(\frac{\omega_m \Gamma(\frac{m}{q}+1)}{\omega_{m+n} \Gamma(\frac{m+n}{q}+1)} \cdot \left((1-t)^mt^n\right)^{\frac{1}{2}-\frac{1}{q}}\right)^\frac{p}{n}.$$
%Since $q\leq 2,$ we choose $t_0=\frac{n}{n+m}$, which gives the best choice for the above expression; thus, we obtain 
\begin{equation}\label{last-ineq}	\int_\Sigma |f|^p\log |f|^p \leq \frac{n}{p}\log\left(\frac{p}{n}\left(\frac{p-1}{e}\right)^{p-1} K_{m,n,q}^\frac{p}{n}
	\int_\Sigma |\nabla^\Sigma f|^p\right),
\end{equation}
where 
\begin{equation}\label{K-sequence}
	K_{m,n,q}=\frac{\omega_m \Gamma(\frac{m}{q}+1)}{\omega_{m+n} \Gamma(\frac{m+n}{q}+1)}\left(\frac{m^mn^n}{(m+n)^{m+n}}\right)^{\frac{1}{2}-\frac{1}{q}}.
\end{equation}
Note that  \eqref{last-ineq} is still a codimension-depending inequality. However, according to Proposition \ref{Alzer-monoton},  the sequence 
$m\mapsto K_{m,n,q} $ is increasing and by the asymptotic property $$\Gamma(r+\alpha)\sim\Gamma(r)r^\alpha\ {\rm as}\ \ r\to \infty$$  with $\alpha\in \mathbb R$, one has that  
\begin{equation}\label{K-limit}
\lim_{m\to \infty}K_{m,n,q}=q^\frac{n}{q}(2\pi)^{-\frac{n}{2}}\left(\frac{e}{n}\right)^{n(\frac{1}{q}-\frac{1}{2})}.
\end{equation}
Combining all these facts, we obtain the codimension-free $L^p$-logarithmic Sobolev  inequality 
$$	\int_\Sigma |f|^p\log |f|^p \leq \frac{n}{p}\log\left( \left(\frac{p^2}{2\pi e n}\right)^{\frac{p}{2}}\int_\Sigma |\nabla^\Sigma f|^p\right),$$
which is precisely inequality \eqref{main-inequality}. 

For $p=2$, inequality \eqref{main-inequality} is sharp; indeed, in this case the constant $\left(\frac{p^2}{2\pi e n}\right)^{\frac{p}{2}}$ is exactly the constant $\mathcal
L_{2,n} = \frac{2}{\pi e n}$ from \eqref{e-sharp-log-Sobolev}. The equality for  $p=2$ is just similar as in Theorem \ref{main-theorem-Brendle}. 
\hfill $\square$

\begin{remark}\rm It is worth to compare the constants in the different versions of the above $L^p$-logarithmic-Sobolev inequalities, namely, \eqref{e-sharp-log-Sobolev}, \eqref{last-ineq}  and \eqref{main-inequality}, respectively. As a result we can conclude that, for every $p\geq 2$,  $n\geq 2$ and $m\geq 1$, we have 
	$$\mathcal L_{p,n}\leq\frac{p}{n}\left(\frac{p-1}{e}\right)^{p-1} K_{m,n,q}^\frac{p}{n}\leq  \left(\frac{p^2}{2\pi e n}\right)^{\frac{p}{2}},$$
	where $q=\frac{p}{p-1}$, whose proof follows by  Proposition \ref{Alzer-monoton} from the Appendix. Moreover, let us note that for $p=2$, the above three constants coincide. 
\end{remark}

%GAMMA(x/q+1)*GAMMA(1/q+1)*GAMMA((x+1)*(1/2)+1)*((x+1)^(x+1)/x^x)^(1/q-1/2)/(GAMMA((1/2)*x+1)*GAMMA(1/2+1)*GAMMA((x+1)/q+1))
%novekvo fuggveny

We conclude this section by stating a codimension-dependent $L^p$-logarithmic Sobolev  inequality ($p\geq 2$) on generic submanifolds: 

\begin{theorem}\label{theorem-general}
	 Let $n\geq2$ and $m\geq 1$ be integers and $p\geq 2$.
	Let $\Sigma$ be a complete $n$-dimensional submanifold of $\mathbb R^{n+m}$.  Then for every function $f\in W_H^{1,p}(\Sigma,d{\rm vol}_\Sigma)$ with $\int_\Sigma |f|^p=1$ we have 
\begin{align*}
\int_\Sigma |f|^p\log |f|^p \leq & \frac{n}{p}\log\left(\left(\frac{2p-3}{p-1}\right)^{\frac{m}{n}(p-1)}A_{n,p} \int_\Sigma |\nabla^\Sigma f|^p\right.\\&\ \ \ \ \ \ \ \ \ \ +\left. \left(\frac{2p-3}{p-1}\right)^{\frac{m+n}{n}(p-1)}B_{n,p}\int_\Sigma |H|^p|f|^p\right),
\end{align*}
where $$A_{n,p}=\left(\frac{p^2}{2\pi e n}\right)^{\frac{p}{2}}\ \ {and}\ \ B_{n,p}=\left(1+n\right)^{\frac{1+n}{n}(\frac{p}{2}-1)} \frac{p-1}{n}\left(\frac{p-1}{pe}\right)^{p-1}\pi^{-\frac{p}{2}}.$$
\end{theorem}

{\it Proof.} 
Since the argument is very similar to the proof of Theorem \ref{main-theorem}, we keep the notations and focus only on the technical differences. 
%We consider the  probability measures $d\mu(x)=|f|^p(x) d{\rm vol}_\Sigma(x)$, $ d\nu_k(y)=C_k^{-1}e^{-\alpha|y|^q}\mathbbm{1}_{B_k}(y)dy,$ and the map $\Phi_k$.  
By the Monge--Amp\`ere equation \eqref{MA-p}, instead of \eqref{log-int}, we obtain an estimate that contains the mean curvature vector $H$, namely, 
 \begin{align}\label{log-int-H-0}
\nonumber	|f|^p(x)\leq&\frac{1}{C_k}\int_{A_x}e^{-\alpha{|\Phi(x,v)|^q}}{\rm det} D\Phi(x,v)dv\\\leq& \frac{1}{M_0}e^{
		-\alpha{t^{1-\frac{q}{2}}}
		\left|\nabla^\Sigma u(x)\right|^q+\Delta_{\rm ac}^\Sigma u(x)-n}\int_{\mathbb R^m}e^{-\alpha{(1-t)^{1-\frac{q}{2}}}|v|^q-\langle H(x),v\rangle}dv.
\end{align}
 Contrary to the arguments from Theorem \ref{main-theorem}, the presence of $H$ prevents the precise computation of the last integral in \eqref{log-int-H-0},  requiring a fine estimate. 
%The latter fact implies that we cannot directly compute the last integral in \eqref{log-int-H} (as in proof of Theorem \ref{main-theorem}), which requires a fine estimate. 
More precisely, if $C_{t,\alpha}:=\alpha(1-t)^{1-\frac{q}{2}}$, it turns out  by Proposition  \ref{elementary-inequality} that for every $x\in \Sigma$ we have
\begin{align*}
-C_{t,\alpha}|v|^q-\langle H(x),v\rangle=&\frac{C_{t,\alpha}}{3-q}\left(-(3-q)|v|^q-\frac{3-q}{C_{t,\alpha}}\langle H(x),v\rangle\right)\\ \leq& \frac{C_{t,\alpha}}{3-q}\left(-|v+w_x|^q+|w_x|^q\right),
\end{align*}
where $w_x\in \mathbb R^m$ solves $\frac{3-q}{C_{t,\alpha}}H(x)=q|w_x|^{q-2}w_x$. Therefore, by \eqref{log-int-H-0} one has for ${\rm vol}_\Sigma$-a.e.\ $x\in \Sigma$ that 
$$|f|^p(x)\leq \frac{\tilde M_{t,\alpha}}{M_0}e^{
	-\alpha{t^{1-\frac{q}{2}}}
	\left|\nabla^\Sigma u(x)\right|^q+\Delta_{\rm ac}^\Sigma u(x)-n+\left(\frac{3-q}{C_{t,\alpha}}\right)^{p-1}\frac{|H(x)|^p}{q^p}},$$
where $$\tilde M_{t,\alpha}=\int_{\mathbb R^m}e^{-\frac{C_{t,\alpha}}{3-q}|v|^q}dv.$$ By  $C_{t,\alpha}=\alpha(1-t)^{1-\frac{q}{2}}$, a similar reasoning as in the proof of Theorem  \ref{main-theorem} gives -- instead of \eqref{later-we-need} -- that
\begin{align}\label{later-we-need-0}	\int_\Sigma |f|^p\log \nonumber|f|^p \leq& \log  \frac{\tilde M_{t,\alpha}}{M_0} -n + \alpha^{1-p}\Biggl[(p-1)^{p-1}t^{(p-1)(\frac{q}{2}-1)} \int_{\Sigma} |\nabla^\Sigma f|^p\\& \ \ \ \ \ \ \ \ \ +\frac{1}{q^p}\left(\frac{3-q}{(1-t)^{1-\frac{q}{2}}}\right)^{p-1}\int_{\Sigma}{|H|^p|f|^p}\Biggl].
\end{align}
%\begin{eqnarray}\label{later-we-need-0}	\int_\Sigma |f|^p\log \nonumber|f|^p &\leq& \log  \frac{\tilde M_{t,\alpha}}{M_0} -n + \alpha^{1-p}\Biggl[(p-1)^{p-1}t^{(p-1)(\frac{q}{2}-1)} \int_{\Sigma} |\nabla^\Sigma f|^p\\&&\ \ \ \ \ \  \ \ \ \ \  \ \  \ \ \ \ \ \ \ \ \  \ \ \ \  +\frac{1}{q^p}\left(\frac{3-q}{(1-t)^{1-\frac{q}{2}}}\right)^{p-1}\int_{\Sigma}{|H|^p|f|^p}\Biggl].
%\end{eqnarray}
Let us use \eqref{Gauss-integral} for $\tilde M_{t,\alpha}$. Minimizing the right hand side w.r.t. $\alpha$ and put $t:=\frac{n}{n+m}$,  we obtain that
$$	\int_\Sigma |f|^p\log |f|^p \leq \frac{n}{p}\log\left( A_{m,n,p}\int_\Sigma |\nabla^\Sigma f|^p+B_{m,n,p}\int_\Sigma |H|^p| f|^p\right),$$
where 
$$A_{m,n,p}=\left(\frac{2p-3}{p-1}\right)^{\frac{m}{n}(p-1)}\frac{p}{n}\left(\frac{p-1}{e}\right)^{p-1} K_{m,n,q}^\frac{p}{n},$$
and 
\begin{align*}
	B_{m,n,p}&=\left(\frac{2p-3}{p-1}\right)^{\frac{m+n}{n}(p-1)}\frac{p-1}{n}\left(\frac{p-1}{pe}\right)^{p-1}\\ &\ \ \ \ \ \ \times \left(\frac{m+n}{m}\right)^{\frac{m+n}{n}(\frac{p}{2}-1)} \left(\frac{\omega_m \Gamma(\frac{m}{q}+1)}{\omega_{m+n} \Gamma(\frac{m+n}{q}+1)}\right)^\frac{p}{n},
\end{align*}
the constant $K_{m,n,q}$ coming from \eqref{K-sequence}. We notice that by the choice of $t=\frac{n}{n+m}$ we minimize the coefficient of the gradient-term, similarly as in the proof of Theorem  \ref{main-theorem}; the joint minimization of the right hand side of \eqref{later-we-need-0}  with respect to the variable $t$ seems to be impossible in the present, general setting, due to its transcendental character. 

 Note that  $\lim_{m\to \infty} A_{m,n,p}=\lim_{m\to \infty} B_{m,n,p}=\infty$ whenever $p\neq 2$. However, by using Propositions \ref{Alzer-monoton} and \ref{Alzer-monoton-2} together with the fact that $x\mapsto (1+\frac{1}{x})^{1+x}$ is decreasing in $(0,\infty)$, it follows that for every $m\geq 1$, $n\geq 2$ and $p\geq 2$, we have
$$A_{m,n,p}\leq \left(\frac{2p-3}{p-1}\right)^{\frac{m}{n}(p-1)}\left(\frac{p^2}{2\pi e n}\right)^{\frac{p}{2}},$$
and 
$$B_{m,n,p}\leq \left(\frac{2p-3}{p-1}\right)^{\frac{m+n}{n}(p-1)}\frac{p-1}{n}\left(\frac{p-1}{pe}\right)^{p-1}\left(1+n\right)^{\frac{1+n}{n}(\frac{p}{2}-1)} \pi^{-\frac{p}{2}},$$
which ends the proof.
\hfill $\square$

\begin{remark}\rm 
If $p=2$, then $\frac{2p-3}{p-1}=1$,  $A_{n,2}=\frac{2}{\pi e n}$ and $B_{n,2}=\frac{1}{2\pi e n}$; in this particular case, Theorem \ref{theorem-general} reduces to the codimension-free Theorem \ref{main-theorem-Brendle}.  

When $p\neq 2$ and $H\neq 0$,   Theorem \ref{theorem-general} is a codimension-dependent $L^p$-logarithmic-Sobolev inequality. In this case,  the sharpness in  Theorem \ref{theorem-general} is a challenging problem, which is similar -- in its nature -- to the sharpness in the isoperimetric inequality on general submanifolds of any codimension, see Brendle \cite{Brendle-1}.
\end{remark}

\section{Hypercontractivity estimates on submanifolds
%	: proof of Theorem \ref{theorem-hypercontractivity}
} \label{section-5}
In this final section we present two applications concerning the hypercontractivity estimate for Hopf--Lax semigroups on submanifolds that are mentioned in the introduction. 

Let us recall that for an $n$-dimensional submanifold $\Sigma$ of $\mathbb R^{n+m}$, a function $u:\Sigma\to \mathbb R$ and for $t>0 $, the  Hopf--Lax semigroup is given by the formula 
%	\begin{equation}\label{inf-convolution-0}
	$$		{\bf Q}_{t}u(x):=\inf_{y\in \Sigma}\left\{u(y)+\frac{d^2_\Sigma(x,y)}{2t}\right\}, \ x\in \Sigma,
	$$
	where $d_\Sigma$ stands for the distance function on $\Sigma.$
	%	Clearly, the admissible $t>0$ and $u:X\to \mathbb R$ are fixed in such a way to ensure that 	${\bf Q}_{t}u(x)\in \mathbb R$;
	It is well known that the Hopf--Lax semigroup gives the solution of the Hamilton--Jacobi equation (\cite{BobkovGL}, \cite{Gentil}, \cite{Villani}) and therefore it is  important to investigate its properties. 
	
	\subsection{Euclidean-type hypercontractivity estimates.}
	 In order to carry out such a study, we need to assume some basic properties of the function $u$. For a similar approach (with a more detailed exposition), see Balogh, Krist\'aly and Tripaldi \cite{BKT}. To be more precise, for some $t_0>0$ and $x_0\in \Sigma$, we introduce the class of functions $u:\Sigma\to \mathbb R$, denoted by $\mathcal  F_{t_0,x_0}(\Sigma)$, which satisfy the following assumptions: 
\begin{itemize}
	\item[(A1)] $u:\Sigma\to \mathbb R$ is smooth and the set $u^{-1}([0,\infty))\subset \Sigma$ is bounded; 
	\item[(A2)]   ${\bf Q}_{t_0}u(x_0)>-\infty;$
	\item[(A3)] there exist $M\in \mathbb R$ and $C_0>2t_0$ such that $$u(x)\geq M-\frac{d_\Sigma^{2}(x,x_0)}{C_0}\ ,\ \ \forall x\in \Sigma.$$
\end{itemize}
Our first result reads as follows: 
\begin{theorem}\label{theorem-hypercontractivity}
	Let $n\geq2$ and $m\geq 1$ be integers and  $\Sigma$ be a complete $n$-dimensional submanifold of $\mathbb R^{n+m}$ with bounded mean curvature vector, i.e. $\|H\|_{L^\infty(\Sigma)}=\sup_{x\in \Sigma}|H(x)|<\infty.$ Let $t_0>0$ and $x_0\in \Sigma$ be fixed. Then for every $0<a<b$, $t\in (0,t_0)$ and $u\in \mathcal F_{t_0,x_0}(\Sigma)$ with $(1+d_\Sigma^{2}(\cdot,x_0))e^{au}\in L^1(\Sigma,d{\rm vol}_\Sigma)$ and $e^\frac{au}{2}\in W^{1,2}(\Sigma,d{\rm vol}_\Sigma)$, one has that 
	\begin{equation}\label{hypercontractivity-estimate}
\|e^{{\bf Q}_tu}\|_{L^b(\Sigma,d{\rm vol}_\Sigma)}\leq \|e^{u}\|_{L^a(\Sigma,d{\rm vol}_\Sigma)}\left(\frac{b-a}{2\pi t}\right)^{\frac{n}{2}\frac{b-a}{ab}}\left(\frac{a}{b}\right)^{\frac{n}{2}\frac{a+b}{ab}}e^{t\frac{\|H\|^2_{L^\infty(\Sigma)}}{6}\frac{a^2+ab+b^2}{a^2b^2}}.
	\end{equation}
Moreover, \eqref{hypercontractivity-estimate} is sharp in the sense that the factor $2\pi$ cannot be replaced by a greater constant, and equality holds for some $t\in (0,t_0)$ and $u\in \mathcal F_{t_0,x_0}(\Sigma)$ with $e^u\in L^a(\Sigma,d{\rm vol}_\Sigma)$ if and only if   $\Sigma$ is isometric to $\mathbb R^n$ and $u(x)=C-\frac{b-a}{2bt}|x-\overline x_0|^2$, $x\in \mathbb R^n,$ for some $C\in \mathbb R$ and $\overline x_0\in \mathbb R^n.$ 
\end{theorem}

	\begin{remark}\rm 
For minimal submanifolds one has $\|H\|_{L^\infty(\Sigma)}=0$ and  \eqref{hypercontractivity-estimate} reduces to 
\begin{equation}\label{hypercontr-estimate-minimal}
\|e^{{\bf Q}_tu}\|_{L^b(\Sigma,d{\rm vol}_\Sigma)}\leq \|e^{u}\|_{L^a(\Sigma,d{\rm vol}_\Sigma)}\left(\frac{b-a}{2\pi t}\right)^{\frac{n}{2}\frac{b-a}{ab}}\left(\frac{a}{b}\right)^{\frac{n}{2}\frac{a+b}{ab}}.
\end{equation}
Thus, for fixed $0<a <b$, the estimate \eqref{hypercontr-estimate-minimal} shows a power decay as $t \to \infty$. 
\end{remark}

%However, in \S\ref{section-3} we state a codimension-depending $L^p$-logarithmic-Sobolev inequality for a general %submanifold in the case when $p\geq 2$, see Theorem \ref{theorem-general}, which reduces to Theorem \ref{main-%theorem-Brendle} for $p=2$ and to Theorem \ref{main-theorem}  for minimal submanifolds. 

%\begin{theorem}
%	Equality holds if and only if  $\Sigma$ is isometric to $\mathbb R^n$ and up to isometry and translation, $f$ is any Gaussian of the form $f_\lambda(x)=(\lambda \pi^{-1})^{\frac{n}{2}} e^{-\lambda |x|^2}$, $x\in \mathbb R^n$, for every $\lambda>0$.
%\end{theorem}
	
%	We apply our logarithmic-Sobolev inequalities to derive Poincar\'e inequalities and hypercontractivity estimates for Hopf--Lax semigroups on submanifolds.

In the sequel, we assume the hypotheses of Theorem \ref{theorem-hypercontractivity} are verified. Namely, let $t_0>0$,  $x_0\in \Sigma$ and  $0<a<b$ fixed, and let $u\in \mathcal F_{t_0,x_0}(\Sigma)$ be such that $(1+d_\Sigma^{2}(\cdot,x_0))e^{au}\in L^1(\Sigma,d{\rm vol}_\Sigma)$ and $e^\frac{au}{2}\in W^{1,2}(\Sigma,d{\rm vol}_\Sigma)$. A similar argument as in \cite{BKT}, essentially based on  Ambrosio, Gigli and Savar\'e \cite{AGS}, shows that the Hopf--Lax semigroup formula
$$		{\bf Q}_{t}u(x):=\inf_{y\in \Sigma}\left\{u(y)+\frac{d^2_\Sigma(x,y)}{2t}\right\}, \ x\in \Sigma,
$$
has the property that $(t,x)\mapsto {\bf Q}_{t}u(x)$ is locally Lipschitz on $(0,t_0)\times \Sigma$ and 
verifies the Hamilton--Jacobi inequality
\begin{equation}\label{HJ-1}
	\frac{d^+}{dt}\textbf{Q}_tu(x)\leq -\frac{|\nabla^\Sigma \textbf{Q}_tu(x)|^2}{2}, \ \ \forall t\in (0,t_0),\ \forall x\in \Sigma, 
\end{equation}
where $\frac{d^+}{dt}$ stands for the right derivative, and for every $x\in \Sigma$ one has
\begin{equation}\label{HJ-2}
	\frac{d^+}{dt}\textbf{Q}_tu(x)\Big|_{t=0}\geq-\frac{|\nabla^\Sigma u(x)|^2}{2}. 
\end{equation}
%where $d_\Sigma$ stands for the distance function on $\Sigma.$
Moreover,  if $\overline t\in (0,t_0)$, for every $c>a$ and any function $q:[0,\overline t]\to [a,c]$ which is of class $C^1$ on $(0,\overline t)$, the function $x\mapsto e^{\frac{q(t)}{2}\textbf{Q}_tu(x)}$ belongs to $W^{1,2}(\Sigma,d{\rm vol}_\Sigma)$  for every $t\in (0,\overline t)$.
Thus, the   function
\begin{align*}
	 F(t):=\bigg(\int_\Sigma e^{q(t){\bf Q}_tu }d{\rm vol}_\Sigma \bigg)^{1/q(t)}=\Vert e^{{\bf Q}_tu}\Vert_{L^{q(t)}(\Sigma,d{\rm vol}_\Sigma)}
\end{align*}
is well-defined and locally Lipschitz on $(0,\overline{t})$. In particular, its right derivative exists at every $t\in (0,\overline t)$ and after an elementary computation it yields that 
\begin{align}\label{right-derivative}
\nonumber	\frac{d^+}{dt} F(t)=&\frac{ F(t)}{q(t)^2}\bigg(\frac{q'(t)}{ F(t)^{q(t)}}\mathcal{E}\big(e^{q(t){\bf Q}_tu }\big)-q'(t)\log\big( F(t)^{q(t)}\big)\\& \ \ \ \ \ \ \ \ \  \ +\frac{q(t)^2}{ F(t)^{q(t)}}\int_\Sigma e^{q(t){\bf Q}_tu }\frac{d^+}{dt}{\bf Q}_tu\bigg),
\end{align}
where 
\begin{align*}
	\mathcal{E}\big(e^{q(t){\bf Q}_tu }\big)=\int_\Sigma e^{q(t){\bf Q}_tu }\log\big(e^{q(t){.\bf Q}_tu }\big)\, d{\rm vol}_\Sigma=q(t)\int_\Sigma e^{q(t){\bf Q}_tu }{\bf Q}_tu\, d{\rm vol}_\Sigma.
\end{align*}

%\textcolor{red}{
	Since  $\|H\|_{L^\infty(\Sigma)}=\sup_{x\in \Sigma}|H(x)|<\infty$,  for every $f\in W^{1,2}(\Sigma,d{\rm vol}_\Sigma)$ with 
$\int_\Sigma f^2=1$ we have 
\begin{equation}\label{inequality-H}
	\int_\Sigma f^2\log f^2 \leq \frac{n}{2}\log\left(\frac{2}{\pi e n }		\left(	\int_\Sigma |\nabla^\Sigma f|^2+\frac{1}{4}\|H\|^2_{L^\infty(\Sigma)}\right)\right),
\end{equation} 
and  \eqref{inequality-H} is sharp, i.e., the constant  $\frac{2}{\pi e n }	$  cannot be improved.
%} 

After this preparatory part, we are ready to provide the proof of Theorem \ref{theorem-hypercontractivity}. \\

{\it Proof of Theorem \ref{theorem-hypercontractivity}.} Let us fix $\overline t\in (0,t_0)$ and consider the function 
	$
	q(t)=\frac{ab}{(a-b)t/\overline{t}+b}, \ t\in (0,\overline t).
$
One clearly has that $q(0)=a$, $q(\overline{t})=b$, and $q$ is of class $C^1$ in $(0,\overline{t})$, with
$
q'(t)=\frac{ab}{[(a-b)t/\overline{t}+b]^2}\cdot\frac{b-a}{\overline{t}}>0\,;
$
thus, $t\mapsto q(t)$ is strictly increasing. 

First of all, according to the the Hamilton--Jacobi inequality \eqref{HJ-1} and relation \eqref{right-derivative}, we obtain for every $t\in (0,\overline t)$ that 
\begin{align}\label{eestimate-derivative}
\nonumber	\frac{\frac{d^+}{dt} F(t)}{F(t)}\le&\frac{q'(t)}{q(t)^2 F(t)^{q(t)}}\bigg(\mathcal{E}\big(e^{q(t){\bf Q}_tu }\big)- F^{q(t)}\log\big( F(t)^{q(t)}\big)\\&\ \ \ \ \ \ \ \ \  \ \ \ \ \ \ \ \ \ \ -\frac{q(t)^2}{q'(t)}\int_\Sigma e^{q(t){\bf Q}_tu }\frac{\vert \nabla^\Sigma{\bf Q}_tu\vert^2}{2}\bigg).
\end{align}

Now, we focus on the terms in the right hand side of \eqref{eestimate-derivative}, where the sharp $L^2$-logarithmic Sobolev inequality \eqref{inequality-H} will be used. 
%Moreover, we introduce the function   $\tilde F:[0,\tilde{t}]\to \mathbb R$ by
%\begin{align*}
%\tilde F(t):=F(t)^{1/q(t)}=\bigg(\int_Xe^{q(t){\bf Q}_tu }d {\sf m}\bigg)^{1/q(t)}=\Vert e^{{\bf Q}_tu}\Vert_{L^{q(t)}(X,{\sf m})}.
%\end{align*}
To do this, for every $t\in (0,\overline t)$ let us define the function 
\begin{align*}
\Sigma \ni	x\mapsto f_t(x):=\frac{e^{\frac{q(t)}{2}{\bf Q}_tu(x)  }}{F(t)^\frac{q(t)}{2}}.
\end{align*}
By the above argument if follows that $f_t$ belongs to $W^{ 1,2}(\Sigma,d{\rm vol}_\Sigma)$  for every 
$t\in (0,\overline t)$ and 
\begin{align*}
	\int_\Sigma f_t^2=1\ ,\ \forall t\in (0,\overline{t}). 
\end{align*}
Due to \eqref{inequality-H}, one has
\begin{equation}\label{LSI-applied-w_t}
	\int_\Sigma f_t^2\log f_t^2 \leq \frac{n}{2}\log\left(\frac{2}{\pi e n }		\left(	\int_\Sigma |\nabla^\Sigma f_t|^2+\frac{1}{4}{\|H\|^2_{L^\infty(\Sigma)}}\right)\right).
\end{equation}
Since we have 
\begin{equation}\label{chain-rule-D}
	\vert \nabla^\Sigma f_t(x)\vert=\frac{e^{{\frac{q(t)}{2}\bf Q}_tu(x) }}{ F(t)^\frac{q(t)}{2}}\cdot\frac{q(t)}{2}\cdot\vert \nabla^\Sigma{\bf Q}_tu(x)\vert\,,
\end{equation}
and 
\begin{equation}\label{identity-transformed}
	\int_\Sigma f_t^2\log f_t^2=\frac{\mathcal{E}\big(e^{q(t){\bf Q}_tu }\big)}{ F(t)^{q(t)}}-\log\big( F(t)^{q(t)}\big),
\end{equation}
it turns out that \eqref{LSI-applied-w_t} is equivalent to 
\begin{align*}
	E:= &\frac{\mathcal{E}\big(e^{q(t){\bf Q}_tu }\big)}{ F(t)^{q(t)}}-\log\big( F(t)^{q(t)}\big)\\  \leq & \frac{n}{2}\log\left(\frac{2}{\pi e n }		\left(	\frac{q^2(t)}{4 F(t)^{q(t)}}\int_\Sigma e^{{{q(t)}\bf Q}_tu(x)} |\nabla^\Sigma {\bf Q}_tu(x)|^2+\frac{1}{4}{\|H\|^2_{L^\infty(\Sigma)}}\right)\right).
\end{align*}
Using the latter estimate and \eqref{eestimate-derivative}, 
we obtain for every $t\in(0,\overline{t})$ that 
\begin{align*}
	\frac{\frac{d^+}{dt} F(t)}{ F(t)}\le&\frac{q'(t)}{q(t)^2}\Bigg[\frac{n}{2}\log\bigg(\frac{1}{2\pi e n }		\Big(	\frac{q^2(t)}{ F(t)^{q(t)}}\int_\Sigma e^{{{q(t)}\bf Q}_tu(x)} |\nabla^\Sigma {\bf Q}_tu(x)|^2\\&\ \ \ \ \ \ \ \  +\|H\|^2_{L^\infty(\Sigma)}\Big)\bigg) -\frac{q(t)^2}{2q'(t) F(t)^{q(t)}}\int_\Sigma e^{q(t){\bf Q}_tu }{\vert \nabla^\Sigma{\bf Q}_tu\vert^2}\Bigg].
\end{align*}
Let us observe that the function 

$$s\mapsto \frac{n}{2}\log\left(\frac{1}{2\pi e n }		\left(	s+{\|H\|^2_{L^\infty(\Sigma)}}\right)\right)-\frac{s}{2q'(t)}, \ \text{defined for} \  s>- {\|H\|^2_{L^\infty(\Sigma)}},$$ has its maximum at the point $s_0=2q'(t)-\|H\|^2_{L^\infty(\Sigma)}$, and its maximal value is equal to 
$$\frac{n}{2}\log\left(\frac{q'(t)}{2\pi e^2}\right)+\frac{\|H\|^2_{L^\infty(\Sigma)}}{2q'(t)}.$$
Therefore, the latter inequality implies that for every $t\in(0,\overline{t})$ we have 
$$\frac{\frac{d^+}{dt}F(t)}{ F(t)}\le \frac{q'(t)}{q(t)^2}\left(\frac{n}{2}\log\left(\frac{q'(t)}{2\pi e^2}\right)+\frac{\|H\|^2_{L^\infty(\Sigma)}}{2q'(t)}\right).$$
Since $t\mapsto  F(t)$ is absolutely continuous and $q(t)=\frac{ab}{(a-b)t/\overline{t}+b}$, an  integration in both sides of the latter inequality over $(0,\overline t)$ implies that 
\begin{align}\label{hyperc-estimate with F}
		F(\overline{t})\le  F(0)\left(\frac{b-a}{2\pi \overline t}\right)^{\frac{n}{2}\frac{b-a}{ab}}\left(\frac{a}{b}\right)^{\frac{n}{2}\frac{a+b}{ab}}e^{\overline t\frac{\|H\|^2_{L^\infty(\Sigma)}}{6}\frac{a^2+ab+b^2}{a^2b^2}}.
\end{align}
Since $\overline t$ was arbitrarily fixed in $(0,t_0)$, the proof of \eqref{hypercontractivity-estimate} is concluded. 

We claim that $2\pi$ cannot be improved in \eqref{hypercontractivity-estimate}; by contradiction, we assume that there exists $C>2\pi $ such that for every $a<b$, $t\in (0,t_0)$ and  $u\in \mathcal F_{t_0,x_0}(\Sigma)$ with $(1+d_\Sigma^{2}(x_0,\cdot))e^{a u}\in  L^1(\Sigma,d{\rm vol}_\Sigma)$ and $e^\frac{a u}{2}\in   W^{ 1,2}(\Sigma,d{\rm vol}_\Sigma)$
  one has
 \begin{equation}\label{hypercontractivity-estimate-C}
 	\|e^{{\bf Q}_tu}\|_{L^b(\Sigma,d{\rm vol}_\Sigma)}\leq \|e^{u}\|_{L^a(\Sigma,d{\rm vol}_\Sigma)}\left(\frac{b-a}{C t}\right)^{\frac{n}{2}\frac{b-a}{ab}}\left(\frac{a}{b}\right)^{\frac{n}{2}\frac{a+b}{ab}}e^{t\frac{\|H\|^2_{L^\infty(\Sigma)}}{6}\frac{a^2+ab+b^2}{a^2b^2}}.
 \end{equation}
Fix a function $u:\Sigma\to \mathbb R$ with the above properties, a number $a>0$, the function $b:=q(t)=a+y t$ with $t\geq 0$, where $y>0$ will be determined later, and 
introduce the function 
$$
\overline F(t):=\bigg(\int_\Sigma e^{q(t){\bf Q}_tu }\, d{\rm vol}_\Sigma\bigg)^{1/q(t)}=\Vert e^{{\bf Q}_tu}\Vert_{L^{q(t)}(\Sigma,d{\rm vol}_\Sigma)},\ t\in [0,t_0).
$$
Having these notations, inequality \eqref{hypercontractivity-estimate-C} reduces  to 
\begin{equation}\label{limit-before}
	\frac{1}{t}{\log \frac{\overline F(t)}{\overline F(0)} }\le \log\bigg[\left(\frac{y}{C }\right)^{\frac{n}{2}\frac{y}{aq(t)}}\left(\frac{a}{q(t)}\right)^{\frac{n}{2}\frac{a+q(t)}{atq(t)}}e^{\frac{\|H\|^2_{L^\infty(\Sigma)}}{6}\frac{a^2+aq(t)+q^2(t)}{a^2q^2(t)}}\bigg],
\end{equation}
for every $ t\in (0,t_0).$  Taking the limit $t\to 0^+$ in \eqref{limit-before}, it follows that 
\begin{align}\label{limit-before-2}
\nonumber	L:=&\frac{\frac{d^+}{dt}\overline F(t)}{\overline F(t)}\bigg|_{t=0}\\ \le & \lim_{t\to 0^+}\log\bigg[\left(\frac{y}{C }\right)^{\frac{n}{2}\frac{y}{aq(t)}}\left(\frac{a}{q(t)}\right)^{\frac{n}{2}\frac{a+q(t)}{atq(t)}}e^{\frac{\|H\|^2_{L^\infty(\Sigma)}}{6}\frac{a^2+aq(t)+q^2(t)}{a^2q^2(t)}}\bigg].
\end{align}
On the one hand, since $$\lim_{t\to 0^+}\left(\frac{a}{q(t)}\right)^{\frac{n}{2}\frac{a+q(t)}{atq(t)}}=e^{-\frac{ny}{a^2}}\ \ {\rm and}\ \ \  \lim_{t\to 0^+}\frac{a^2+aq(t)+q^2(t)}{a^2q^2(t)}=\frac{3}{a^2},$$
the estimate \eqref{limit-before-2} reduces to 
\begin{equation}\label{limit-before-2-bis}
	L=\frac{\frac{d^+}{dt}\overline F(t)}{\overline F(t)}\bigg|_{t=0}\le \frac{ny}{2a^2}\log\left(\frac{y}{c}\right)-\frac{ny}{a^2}+\frac{\|H\|^2_{L^\infty(\Sigma)}}{2a^2}.
\end{equation}
On the other hand, by \eqref{right-derivative} and \eqref{HJ-2} one has that 
\begin{align*}
	\nonumber L=	\frac{\frac{d^+}{dt}\overline F(t)}{\overline F(t)}\bigg|_{t=0}=&\frac{1}{q(t)^2}\bigg(\frac{q'(t)}{ \overline F(t)^{q(t)}}\mathcal{E}\big(e^{q(t){\bf Q}_tu }\big)-q'(t)\log\big( \overline F(t)^{q(t)}\big)\\&\ \ \ \ \ \  \ \ +\frac{q(t)^2}{\overline F(t)^{q(t)}}\int_\Sigma e^{q(t){\bf Q}_tu }\frac{d^+}{dt}{\bf Q}_tu\bigg)\Bigg|_{t=0}\\ \geq& \nonumber \frac{y}{a^2\|e^u\|^a_{L^{a}(\Sigma)}}\bigg (\mathcal{E}\big(e^{a u}\big)-\|e^u\|^a_{L^{a}(\Sigma)}\log (\|e^u\|^a_{L^{a}(\Sigma)})\\&\ \ \ \ \ \ \ \ \ \ \ \ \ \ \ \ \ \ \ -\frac{a^2}{2y}\int_\Sigma e^{a u}|\nabla^\Sigma u|^2\bigg).
\end{align*}
By assumption, $e^\frac{a u}{2}\in   W^{ 1,2}(\Sigma,d{\rm vol}_\Sigma)$, therefore  
$$f:=\frac{e^\frac{a u}{2}}{\|e^u\|^\frac{a}{2}_{L^{a}(\Sigma)}}\in  W^{ 1,2}(\Sigma,d{\rm vol}_\Sigma)\ \ {\rm and}\ \ \displaystyle\int_\Sigma f^{2}=1. $$
Using this notation and combining the latter estimate,  inequality   \eqref{limit-before-2-bis} implies that 
$$\int_\Sigma f^2\log f^2\leq \frac{2}{y}\int_\Sigma |\nabla^\Sigma f|^2+\frac{\|H\|^2_{L^\infty(\Sigma)}}{2y}+\frac{n}{2}\log\left(\frac{y}{C}\right)-n.$$
Due to the fact that this inequality is valid for every $y>0$, we can minimize the right hand side with respect to $y$, obtaining the minimum point $y_0:=\frac{4}{n}\left(	\int_\Sigma |\nabla^\Sigma f|^2+\frac{1}{4}\|H\|^2_{L^\infty(\Sigma)}\right);$ accordingly, it follows that
$$\int_\Sigma f^2\log f^2\leq \frac{n}{2}\log\left(\frac{4}{Cen }		\left(	\int_\Sigma |\nabla^\Sigma f|^2+\frac{1}{4}\|H\|^2_{L^\infty(\Sigma)}\right)\right).$$
Since $C>2\pi$, the latter inequality contradicts the sharpness of the constant $\frac{2}{\pi e n }	$ in \eqref{inequality-H}. 

Let us assume the equality holds in \eqref{hypercontractivity-estimate} for  some $t\in (0,t_0)$ and $u\in \mathcal F_{t_0,x_0}(\Sigma)$ with $e^u\in L^a(\Sigma,d{\rm vol}_\Sigma)$. Consequently, we should have equality in the $L^2$-logarithmic Sobolev inequality \eqref{LSI-applied-w_t} for $f_t$ as well; in particular, by the equality case in Theorem \ref{main-theorem-Brendle}  it follows that $\Sigma$ is isometric to the Euclidean space $\mathbb R^n$. 
% and for some $\alpha>0$ and $x_0\in\mathbb R^n,$  one has $f_t(x)=(\frac{\alpha}{\pi})^{\frac{n}{4}} e^{-\frac{\alpha}{2} |x-x_0|^2}$, $x\in \mathbb R^n$. 
At this point we can use the recent result by Balogh, Don and Krist\'aly \cite[Theorem 1.3]{BDK}, which characterizes the equality case in the hypercontractivity estimate in Euclidean spaces, -- in particular, in \eqref{hypercontr-estimate-minimal} -- obtaining that $u(x)=C-\frac{b-a}{2bt}|x-\overline x_0|^2$, $x\in \mathbb R^n,$ for some $C\in \mathbb R$ and $\overline x_0\in \mathbb R^n.$ 
  \hfill $\square$
  
  \subsection{Gaussian-type hypercontractivity estimates.}
We now focus on a sharp hypercontractivity estimate with respect to the Gaussian measure
$$d\gamma(x)=\left(4\pi\right)^{-\frac{n}{2}}e^{-\frac{|x|^2}{4}}d{\rm vol}_\Sigma(x),\ \ x\in \Sigma,$$
whenever $\Sigma\subset \mathbb R^{n+m}$ is a\textit{ self-similar shrinker}, i.e.,   $H+\frac{x^\perp}{2}=0$ holds on $\Sigma.$ Before stating our result we recall two essential properties of 
self-similar shrinkers, which are crucial in our argument. 

Let us note first that  since $H+\frac{x^\perp}{2}=0$  on $\Sigma,$  by Corollary \ref{corollary-gauss}  it follows that for every
$\varphi\in W_{H,1}^{1,2}(\Sigma,d\gamma)=\{ \varphi\in W^{1,2}(\Sigma,d\gamma):|x| \varphi\in L^2(\Sigma,d\gamma)\}$ with $\int_\Sigma \varphi^2d\gamma=1$ one has
\begin{equation}\label{parametric-gauss-1}
	\int_\Sigma \varphi^2\log\varphi^2 d\gamma\leq 4
	\int_\Sigma |\nabla^\Sigma \varphi|^2d\gamma.
\end{equation}
In addition, a similar argument as in \cite{BKT} shows that the constant $4$ is sharp in \eqref{parametric-gauss-1}.

Furthermore, due to Ding and Xin \cite[Theorem 1.1]{DX},  the surface measure of any $n$-dimensional self-similar shrinker  $\Sigma$  of $\mathbb R^{n+m}$ has \textit{Euclidean volume growth}: there exists a constant $C(n)>0$ depending only on $n$ such that 
	\begin{equation}\label{growth-estimate}
		{\rm vol}_\Sigma(B_r\cap \Sigma)\leq C(n) r^n, \ \ \forall r\geq 1,
	\end{equation}
	 where $B_r$ is the usual ball in $\mathbb R^{n+m}$ 
with radius $r$ and centered at the origin. A useful consequence of \eqref{growth-estimate} is that 
\begin{equation}\label{growth-beta-self-shrinker}
	|x|^\beta e^{C|x|^\theta}\in L^1(\Sigma,d\gamma),\ \ \ \forall \beta,C\geq 0,\ \theta\in (0,2).
\end{equation}
Indeed, without loss of generality, after suitable translation, we may assume that $B_r\cap \Sigma =\emptyset$ for sufficiently large $r>0$; thus, by using the layer cake representation, one has by \eqref{growth-estimate} that  
 \begin{align*}
 	\int_{\Sigma} |x|^\beta e^{C|x|^\theta} d\gamma =&\left(4\pi\right)^{-\frac{n}{2}} \int_{\Sigma}|x|^\beta  e^{C|x|^\theta-\frac{|x|^2}{4}} d{\rm vol}_{\Sigma} \leq C' \int_{\Sigma}|x|^\beta e^{-\frac{|x|^2}{5}}d{\rm vol}_{\Sigma}\\ \leq&C'' \int_0^{\infty}{\rm vol}_\Sigma(B_r\cap \Sigma) (r^{\beta-1}+r^{\beta+1}) e^{-\frac{r^2}{5}}  dr\\  \leq& C'' C(n) \int_0^{\infty} (r^{n+\beta-1}+r^{n+\beta+1}) e^{-\frac{r^2}{5}}  dr  <+\infty,
 \end{align*}
for some $C',C''>0$. 
%Now, using  of the integral on the right side of the above inequality, combined with  we can further estimate 
%\begin{equation} \label{start-integral0} 
%	\int_{\Sigma} e^{au} d\gamma \leq  \leq C'' \int_0^{\infty} r^{n+1} e^{-\frac{r^2}{5}} dr < \infty,  
%\end{equation} 

Armed with this preparatory part, our sharp hypercontractivity estimate on self-similar shrinkers reads as follows. 
\begin{theorem}\label{them-hyper-Gauss}
Let	$\Sigma\subset \mathbb R^{n+m}$ be an $n$-dimensional self-similar shrinker.  Assume that 
$u: \Sigma \to \mathbb R$ is a smooth function with the property that there exist $C_1, C_2 >0$ and $0< \theta < 2$ such that 
\begin{equation}  \label{growth-cond}  |u(x)| \leq C_1 + C_2 |x|^{\theta}, \ x \in \Sigma.
\end{equation} 
Then for every $a>0$ and $t>0$, one has 
	\begin{equation}\label{hypercontractivity-estimate-Gauss-1}
		\|e^{{\bf Q}_tu}\|_{L^{a+\frac{t}{2}}(\Sigma,d\gamma)}\leq \|e^{u}\|_{L^a(\Sigma,d\gamma)}.
	\end{equation}
	Moreover, the factor $2$ in the $L^{a+\frac{t}{2}}(\Sigma,d\gamma)$-norm cannot be replaced by a smaller constant, and equality holds in \eqref{hypercontractivity-estimate-Gauss-1} for some $a>0$, $t>0$ and $u:\Sigma \to \mathbb R$ satisfying \eqref{growth-cond} if and only if $\Sigma$ is isometric to $\mathbb R^n$ and $u(x) = \langle x, x_0 \rangle + C_0$ for some $x_0 \in \mathbb R^n$ and $C_0 \in \mathbb R.$ 
\end{theorem}

{\it Proof.} The line of the proof is similar to the one presented by Bobkov, Gentil and Ledoux \cite[Theorem 2.1]{BobkovGL}. In order to carry out the formal calculations as in \cite{BobkovGL} we need to establish the necessary regularity properties.

Fix $a>0$ and pick a smooth function $u:\Sigma \to \mathbb R$ satisfying the growth condition \eqref{growth-cond}. According to \eqref{growth-cond} and \eqref{growth-beta-self-shrinker}, we observe that $e^u \in {L^a(\Sigma,d\gamma)}$ and in particular the left side of \eqref{hypercontractivity-estimate-Gauss-1} is finite. 
%	To see this, let us use the growth condition \eqref{growth-cond} to obtain
%	$$ \int_{\Sigma} e^{au} d\gamma \leq C \int_{\Sigma} e^{aC_2 |x|^{\alpha}-\frac{|x|^2}{4}} d{\rm vol}_{\Sigma} \leq \hat{C} \int_{\Sigma} e^{-\frac{|x|^2}{5}}d{\rm vol}_{\Sigma},$$
%for some constants $C, \hat{C}>0$ depending only on the values of $C_1, C_2, a>0$. 
%Now, using the layer cake representation of the integral on the right side of the above inequality, combined with \eqref{growth-estimate} we can further estimate 
%\begin{equation} \label{start-integral} 
%\int_{\Sigma} e^{au} d\gamma \leq C' \int_0^{\infty}{\rm vol}_\Sigma(B_r\cap \Sigma) e^{-\frac{r^2}{5}} r dr \leq C'' \int_0^{\infty} r^{n+1} e^{-\frac{r^2}{5}} dr < \infty,  
%\end{equation}
Furthermore, we notice that by the definition of ${\bf Q}_tu$ we have ${\bf Q}_tu(x) \leq u(x)$ for every $t>0$ and $x\in \Sigma$;  thus we can apply a similar argument to conclude that $ e^{{\bf Q}_tu} \in L^{a+\frac{t}{2}}(\Sigma,d\gamma)$.  
For every $t>0$, define the functions $q(t)=a+\frac{t}{2}$  and 
\begin{align*}
	F(t)=\bigg(\int_\Sigma e^{q(t){\bf Q}_tu }d\gamma \bigg)^{1/q(t)}=\Vert e^{{\bf Q}_tu}\Vert_{L^{q(t)}(\Sigma,d\gamma)}.
\end{align*}
The above consideration just showed that $F$ is well-defined for all $t\geq 0$. Moreover, a slight modification of the arguments in \cite[Proposition 4.1]{BKT} (see also \cite{BobkovGL}) shows that $(t,x)\mapsto {\bf Q}_tu(x)$ is locally Lipschitz on $(0,\infty)\times \Sigma$, and relations \eqref{HJ-1} and \eqref{HJ-2} hold for every $t>0$ and $x\in \Sigma$.
This fact also implies  that $F^q$, and $F$, are locally Lipschitz on $(0,\infty)$. In particular, $-\infty<{\frac{d^+}{dt} F(t)^{q(t)}}$ for every $t>0$; thus, by the Hamilton--Jacobi inequality \eqref{HJ-1} one has that $$-\infty<\int_\Sigma e^{q(t){\bf Q}_tu } \frac{d^+}{dt}{\bf Q}_tu  d\gamma \leq -\frac{1}{2}\int_\Sigma e^{q(t){\bf Q}_tu } |\nabla^\Sigma \textbf{Q}_tu(x)|^2  d\gamma.$$
The latter relation together with  $e^{q(t){\bf Q}_tu }\in L^1(\Sigma,d\gamma)$ implies that $x\mapsto e^{\frac{q(t)}{2}{\bf Q}_tu(x)  }$ belongs to $W^{1,2}(\Sigma,d\gamma)$ for every $t>0$. 
{Moreover, using \eqref{growth-beta-self-shrinker} and \eqref{growth-cond} as before, it follows that  $x\mapsto |x| e^{\frac{q(t)}{2}{\bf Q}_tu(x)  }\in L^2(\Sigma,d\gamma)$ for every $t>0$. 
%Indeed, by \eqref{growth-estimate} and $B_2\cap \Sigma =\emptyset$, a similar arguments as in \eqref{growth-est} shows that 
%\begin{eqnarray*}\label{weighted-Sobolev-ben-van}
	%\int_\Sigma |x|^2e^{q(t){\bf Q}_tu(x) }d\gamma(x) &\leq& e^{Mq(t) }\int_\Sigma |x|^2d\gamma(x)=e^{Mq(t) }\left(4\pi\right)^{-\frac{n}{2}}\int_\Sigma |x|^2 e^{-\frac{|x|^2}{4}}d{\rm vol}_\Sigma(x)\\&\leq & \frac{1}{2}e^{Mq(t) }\left(4\pi\right)^{-\frac{n}{2}}\int_0^\infty {\rm vol}_\Sigma(B_r\cap \Sigma) r(r^2-4) e^{-\frac{r^2}{4}}dr \\&<&2(n+2)e^{Mq(t)} \frac{C(n)}{\omega_n} <\infty.
%\end{eqnarray*}
%}
Therefore, if 
$
\Sigma \ni	x\mapsto \varphi_t(x):=\frac{e^{\frac{q(t)}{2}{\bf Q}_tu(x)  }}{F(t)^\frac{q(t)}{2}},
$
it turns out that $\varphi_t\in W^{1,2}_{H,1}(\Sigma,d\gamma)$ and $\int_\Sigma \varphi_t^2d\gamma =1$; in particular, $\varphi_t$ verifies \eqref{parametric-gauss-1}. After a similar computation as in the proof of Theorem \ref{theorem-hypercontractivity} we have that
\begin{equation}\label{monotonicity}
	\frac{\frac{d^+}{dt} F(t)}{F(t)}\le\frac{q'(t)}{q(t)^2 }\bigg(\int_\Sigma \varphi_t^2\log\varphi_t^2d\gamma-4\int_\Sigma 	\vert \nabla^\Sigma \varphi_t\vert^2 d\gamma\bigg)\leq 0.
\end{equation}
This relation implies that $t\mapsto F(t)$ is non-increasing on $(0,\infty)$; in particular, $F(t)\leq F(0)$ for every $t>0$, which is precisely inequality \eqref{hypercontractivity-estimate-Gauss-1}. 

%
% \begin{align}\label{eestimate-derivative-1}
% 	\frac{\frac{d^+}{dt} F(t)}{F(t)}\le\frac{q'(t)}{q(t)^2 F(t)^{q(t)}}\bigg(\mathcal{E}_\gamma\big(e^{q(t){\bf Q}_tu }\big)- F^{q(t)}\log\big( F(t)^{q(t)}\big)-\frac{q(t)^2}{q'(t)}\int_\Sigma e^{q(t){\bf Q}_tu }\frac{\vert \nabla^\Sigma{\bf Q}_tu\vert^2}{2}d\gamma\bigg),
% \end{align}
%where 
%\begin{align*}
%	\mathcal{E}_\gamma\big(e^{q(t){\bf Q}_tu }\big)=\int_\Sigma e^{q(t){\bf Q}_tu }\log\big(e^{q(t){.\bf Q}_tu }\big)d\gamma=q(t)\int_\Sigma e^{q(t){\bf Q}_tu }{\bf Q}_tud\gamma.
%\end{align*}

%
%
% moreover,   
%$$\frac{\mathcal{E}_\gamma\big(e^{q(t){\bf Q}_tu }\big)}{F^{q(t)}}- \log\big( F(t)^{q(t)}\big)=	\int_\Sigma \varphi_t^2\log\varphi_t^2d\gamma$$
%and $$
%	\vert \nabla^\Sigma \varphi_t(x)\vert=\frac{e^{{\frac{q(t)}{2}\bf Q}_tu(x) }}{ F(t)^\frac{q(t)}{2}}\cdot\frac{q(t)}{2}\cdot\vert \nabla^\Sigma{\bf Q}_tu(x)\vert.
%$$
%Thus, by \eqref{eestimate-derivative-1} and \eqref{parametric-gauss-1} one has that
%$$\frac{\frac{d^+}{dt} F(t)}{F(t)}\le\frac{q'(t)}{q(t)^2 }\bigg(\int_\Sigma \varphi_t^2\log\varphi_t^2d\gamma-4\int_\Sigma 	\vert \nabla^\Sigma \varphi_t\vert^2 d\gamma\bigg)\leq 0.$$
%Thus $t\mapsto F(t)$ is non-increasing; in particular, $F(t)\leq F(0)$ for every $t>0$, which is precisely inequality \eqref{hypercontractivity-estimate-Gauss-1}. 

	Let us assume that one can replace the number 2 in  \eqref{hypercontractivity-estimate-Gauss-1} by a smaller constant, i.e., there exists $b<2$ such that for every $a,t>0$ and smooth, bounded function $u:\Sigma\to \mathbb R$, one has
\begin{equation}\label{sharpness-question}
	\|e^{{\bf Q}_tu}\|_{L^{a+\frac{t}{b}}(\Sigma,d\gamma)}\leq \|e^{u}\|_{L^a(\Sigma,d\gamma)}.
\end{equation}
First, since $e^\frac{a u}{2}\in   W_{H,1}^{ 1,2}(\Sigma,d\gamma)$, one has that
$$\varphi:=\frac{e^\frac{a u}{2}}{\|e^u\|^\frac{a}{2}_{L^{a}(\Sigma,d\gamma)}}\in  W_{H,1}^{ 1,2}(\Sigma,d\gamma)\ \ {\rm and}\ \ \displaystyle\int_\Sigma \varphi^{2}d\gamma=1. $$
Moreover, if 
\begin{align*}
	\overline F(t)=\bigg(\int_\Sigma e^{q(t){\bf Q}_tu }d\gamma \bigg)^{1/q(t)}=\Vert e^{{\bf Q}_tu}\Vert_{L^{q(t)}(\Sigma,d\gamma)},
\end{align*}
where $q(t)=a+\frac{t}{b}$, the inequality \eqref{sharpness-question} means that  $\overline F(t)\leq \overline F(0)$ for every $t>0$; in particular, $\frac{d^+}{dt} \overline F(t)\big|_{t=0}\leq 0.$ This relation together with \eqref{HJ-2}  yields 
%\begin{eqnarray}
%	\nonumber 0&\geq&	\frac{\frac{d^+}{dt}\overline F(t)}{\overline F(t)}\bigg|_{t=0}
%%	\\&=&\nonumber \frac{1}{q(t)^2}\bigg(\frac{q'(t)}{ \overline F(t)^{q(t)}}\mathcal{E}_\gamma\big(e^{q(t){\bf Q}_tu }\big)-q'(t)\log\big( \overline F(t)^{q(t)}\big)\\&&\ \ \ \ \ \ \ \  \ +\nonumber\frac{q(t)^2}{\overline F(t)^{q(t)}}\int_\Sigma e^{q(t){\bf Q}_tu }\frac{d^+}{dt}{\bf Q}_tu d\gamma\bigg)\Bigg|_{t=0}
%	\\&\geq &\nonumber \frac{1}{a^2b\|e^u\|^a_{L^{a}(\Sigma,d\gamma)}}\bigg(\mathcal{E}_\gamma\big(e^{a u}\big)-\|e^u\|^a_{L^{a}(\Sigma,d\gamma)}\log (\|e^u\|^a_{L^{a}(\Sigma,d\gamma)})\\&&\nonumber \ \ \ \ \ \ \  \  \ \ \ \ \ \ \  \ \ \ \ \ \ \ \  \ -\frac{a^2b}{2}\int_\Sigma e^{a u}|\nabla^\Sigma u|^2d\gamma\bigg).
%\end{eqnarray}
%Since $e^\frac{a u}{2}\in   W^{ 1,p}(\Sigma,d\gamma)$, one has that
%$$\varphi:=\frac{e^\frac{a u}{2}}{\|e^u\|^\frac{a}{2}_{L^{a}(\Sigma,d\gamma)}}\in  W^{ 1,2}(\Sigma,d\gamma)\ \ {\rm and}\ \ \displaystyle\int_\Sigma \varphi^{2}d\gamma=1. $$
%The above estimate can be equivalently written as 
$$\int_\Sigma \varphi^2\log\varphi^2d\gamma-2b\int_\Sigma 	\vert \nabla^\Sigma \varphi\vert^2 d\gamma\leq 0.$$
Since $b<2$, the latter inequality contradicts the sharpness of the constant $4$ in \eqref{parametric-gauss-1}. 

It remains to characterize the equality case; thus, we assume that equality holds in \eqref{hypercontractivity-estimate-Gauss-1} for some $a>0$, $t>0$ and smooth function $u:\Sigma \to \mathbb R$. In particular, it follows that $F(s)= F(0)$ for every $s\in (0,t],$ thus by \eqref{monotonicity} one has 
$$\int_\Sigma \varphi_s^2\log\varphi_s^2d\gamma-4\int_\Sigma 	\vert \nabla^\Sigma \varphi_s\vert^2 d\gamma=0, \ \ \forall s\in (0,t].$$
  By the equality case in Corollary \ref{corollary-gauss}, it turns out that $\Sigma$ is isometric to the Euclidean space $\mathbb R^n$ and  $\varphi_s(x)=e^{\frac{1}{4} \langle x,x_s\rangle -\frac{1}{8}|x_s|^2}$, $x\in \mathbb R^n$, for some $x_s\in \mathbb R^n.$ Accordingly, it follows that
 $$\frac{e^{\frac{q(s)}{2}{\bf Q}_su(x)  }}{F(s)^\frac{q(s)}{2}}=e^{\frac{1}{4} \langle x,x_s\rangle -\frac{1}{8}|x_s|^2},\ \ \forall s\in (0,t],\ x\in \mathbb R^n.$$
 
Changing notation, the above relation implies that  for every $s\in (0,t]$ one has
 ${\bf Q}_su(x)= \langle x, x_s \rangle + c_s$ for some $x_s \in \mathbb{R}^n$ and $c_s \in \mathbb{R}$. Let us use the semigroup property of the inf-convolution $ s \mapsto {\bf Q}_s$, i.e.,  for $t>s\geq 0$ and all $x \in \Sigma,$ we have the relation 
 \begin{equation} \label{semigroup} 
  {\bf Q}_{t-s}{\bf Q}_s u(x) = {\bf Q}_tu(x).
  \end{equation} 
 Using the defining formula of ${\bf Q}_{t-s}$ and the expression ${\bf Q}_su(x)= \langle x, x_s \rangle + c_s,$ we can calculate explicitly the left side of the above relation, namely  $$ {\bf Q}_{t-s}{\bf Q}_s u(x) = \langle x, x_s \rangle + c_s - |x_s|^2 \frac{t-s}{2}.$$ 
 This relation with \eqref{semigroup} implies
 $$ \langle x, x_s \rangle + c_s -|x_s|^2 \frac{t-s}{2} = \langle x, x_t \rangle + c_t.$$
 In particular, we immediately obtain that $x_s = x_t$, i.e., $x_t= x_0$ (it does not depend on $t$) and thus 
 ${\bf Q}_s u(x) = \langle x, x_0\rangle + c_s$, $s\in (0,t]$. In addition, the latter relation forces
 $$  c_s +\frac{s}{2}|x_0|^2  =   c_t+\frac{t}{2}|x_0|^2,\ s\in (0,t] .$$
 Therefore, $c_s +\frac{s}{2}|x_0|^2\equiv C_0$ for some $C_0\in \mathbb R$ and for every  $s\in (0,t]$.   It is easy to see that   $\lim_{s \to 0^+}{\bf Q}_s u(x) = u(x)$; thus   $\lim_{s \to 0^{+}} c_s =C_0 $ and  
 $u(x) = \langle x, x_0\rangle + C_0$ as claimed.  
%By the boundedness of $(s,x)\mapsto {\bf Q}_su(x)$, see \eqref{u-bounded},  the latter relation necessarily implies that $x_s=0$ for every $s\in (0,t]$. Therefore, for every $s\in (0,t]$ there exists $C_s\in \mathbb R$ such that ${\bf Q}_su(x)=C_s$ for every $x\in \mathbb R^n$. Since $C_s={\bf Q}_su\leq u$, it follows that
%$$e^{C_s}\leq \|e^{u}\|_{L^a(\mathbb R^n,d\gamma)}=F(0)=F(s)=\|e^{{\bf Q}_su}\|_{L^{a+\frac{s}%{2}}(\mathbb R^n,d\gamma)}= \|e^{C_s}\|_{L^{a+\frac{s}{2}}(\mathbb R^n,d\gamma)}=e^{C_s};$$  
%in particular, we have for every $s\in (0,t]$ that $u(x)=C_s$ for $\gamma$-a.e. $\mathbb R^n$. By the smoothness of $u$, we have in fact that $u\equiv C_0$ for some $C_0\in \mathbb R$,  
%which completes the proof.
\hfill $\square$

\begin{remark}\label{remark-theta}\rm Let us comment on the range of the parameter $\theta$ in the growth condition \eqref{growth-cond}.
	
	On the one hand,  as the proof of Theorem \ref{them-hyper-Gauss} shows, assumption $\theta\in (0,2)$ provides an appropriate range to prove the  hypercontractivity estimate \eqref{hypercontractivity-estimate-Gauss-1}. In addition, if $\theta<1$ in \eqref{growth-cond} (i.e., $u$ is sublinear), equality holds in \eqref{hypercontractivity-estimate-Gauss-1} if and only if $\Sigma$ is isometric to $\mathbb R^n$ and $u$ is constant;  indeed, by  $u(x) = \langle x, x_0\rangle + C_0$ with $x_0\neq 0$, the growth condition \eqref{growth-cond} fails for  $\theta<1$. 
	
	On the other hand, the case $\theta\geq 2$ is much more  sensitive. Indeed,  for simplicity, if we assume that $\Sigma=\mathbb R^n$ and $u(x)=-|x|^\theta$ with $\theta>2$,  it turns out that ${\bf Q}_t u(x)=-\infty$ for every $t>0$ and $x\in \Sigma$, which provides no reasonable hypercontractivity estimates.  However, the threshold value $\theta=2$ in \eqref{growth-cond} still guarantees the validity of \eqref{hypercontractivity-estimate-Gauss-1}, whenever $a,t,C_2>0$ verify $(a+\frac{t}{2})C_2<\frac{1}{4}$. Indeed, in this range of parameters, we have ${\bf Q}_t u(x)\in \mathbb R$ for every $x\in \Sigma$, while $e^u \in {L^a(\Sigma,d\gamma)}$ and $|x|^\beta e^{(a+\frac{t}{2}){\bf Q}_tu}\in L^1(\Sigma,d\gamma)$ for every $\beta\geq 0;$ the rest of the proof is similar as above. 
\end{remark}

\appendix
\section{}

In this final part we collect a variety of technical results that could not be found in the literature, but they are used in the proofs of our main results.  These include  monotonicity properties of the Gamma function, as well as a non-symmetric parallelogram-type inequality.

		\subsection{Monotonicity of certain special functions related to the $\Gamma$ function}
		
	The following results are used in the proof of Theorem \ref{main-theorem}. They express certain monotonicity properties of special functions related to the $\Gamma$ function.

First, we recall that if $k\geq 1$ is an integer and $q,\alpha>0$,  a standard computation shows that
	\begin{equation}\label{Gauss-integral}
		C_{\alpha,k}:=\int_{\mathbb R^k}e^{-\alpha |x|^q}dx=\alpha^{-\frac{k}{q}}\omega_k\Gamma\left(\frac{k}{q}+1\right).
	\end{equation}	
	We first focus on the monotonicity of the sequence 
	$m\mapsto K_{m,n,q} $, appearing in \eqref{K-sequence}. By using the expression $\omega_k= \frac{\pi^{\frac{k}{2}}}{\Gamma(\frac{k}{2}+1)} $ we can write $K_{m,n,q}$ only by means of the $\Gamma$ function, namely, 
	$$	K_{m,n,q}=\frac{1}{\pi^\frac{n}{2}}\frac{ \Gamma(\frac{m+n}{2})}{ \Gamma(\frac{m}{2})}\frac{ \Gamma(\frac{m}{q})}{ \Gamma(\frac{m+n}{q})}\left(\frac{m^mn^n}{(m+n)^{m+n}}\right)^{\frac{1}{2}-\frac{1}{q}},$$
	where we used the recurrence relation $\Gamma(x+1)=x\Gamma(x),$ $x>0.$
	  Accordingly, we see that the desired monotonicity of $m\mapsto K_{m,n,q} $ follows from the following result: 
	
	\begin{proposition}\label{Alzer-monoton}
		For every $t\geq 0$ and $q\in (1,2]$, the function 
		$$x\mapsto \frac{ \Gamma(\frac{x+t}{2})}{ \Gamma(\frac{x+t}{q})}\frac{ \Gamma(\frac{x}{q})}{ \Gamma(\frac{x}{2})}\left(\frac{x^x}{(x+t)^{x+t}}\right)^{\frac{1}{2}-\frac{1}{q}}$$ is increasing on $(0,\infty)$. 
	\end{proposition}

{\it Proof.} Clearly, for $t=0$ or $q=2$ the statement trivially holds. For the case of general $t, q$ we shall prove that the logarithmic derivative of the above function is non-negative. If $\psi=\frac{\Gamma'}{\Gamma}$ denotes the digamma function, we should prove that for every $x>0$, $t\geq 0$ and $q\in (1,2]$ one has
\begin{align}\label{monotonic}
\nonumber	h(x,t,q):=&\frac{1}{2}\psi\left(\frac{x+t}{2}\right)-\frac{1}{q}\psi\left(\frac{x+t}{q}\right)+\frac{1}{q}\psi\left(\frac{x}{q}\right)\\&-\frac{1}{2}\psi\left(\frac{x}{2}\right)+ \left(\frac{1}{2}-\frac{1}{q}\right)\log\frac{x}{x+t}\\ \geq& 0.\nonumber
\end{align}
On the one hand, we recall the result of Alzer \cite[Theorem 4]{Alzer}, which asserts that
$$x(x\psi(x))''<1,\ \forall x>0.$$
This inequality implies  that the function $x\mapsto x^2 \psi'(x)-x$ is decreasing on $(0,\infty)$. Applying this monotonicity property for the points $\frac{x+t}{q}\geq \frac{x+t}{2}$  for  $x>0$, $t\geq 0$ and $q\in (1,2]$, it follows that
$$\frac{1}{4}\psi'\left(\frac{x+t}{2}\right)-\frac{1}{q^2}\psi'\left(\frac{x+t}{q}\right)-\left(\frac{1}{2}-\frac{1}{q}\right)\frac{1}{x+t}\geq 0.$$
On the other hand, observe that $\frac{\partial}{\partial t}h(x,t,q)$ is exactly the left hand side of the latter expression; thus,  $\frac{\partial}{\partial t}h(x,t,q)\geq 0.$ In particular, $t\mapsto 	h(x,t,q)$ is increasing. Thus, for every $x>0$, $t\geq 0$ and $q\in (1,2]$ we obtain $h(x,t,q)\geq h(x,0,q)=0$, which concludes the proof of \eqref{monotonic}. 
\hfill $\square$\\

Another important monotonicity result reads as follows: 

\begin{proposition}\label{Alzer-monoton-2}
	For every  $r\geq 1$, the function 
	$x\mapsto \frac{ \Gamma(rx)}{ \Gamma(x)}$ is increasing on $[\frac{1}{2},\infty)$. 
\end{proposition}

{\it Proof.} Denoting again by $\psi = \frac{\Gamma'}{\Gamma}$ the digamma function, let us note that the statement of the proposition is equivalent to the inequality $r\psi(rx) - \psi(x) \geq 0 $ for $r\geq 1$ and $x\geq \frac{1}{2} $. 

To check the above inequality we start from the result of Alzer \cite[Theorem 4]{Alzer}, which states  that  
$(x\psi(x))''>0$ for every $ x>0.$
This implies that $x\mapsto \psi(x)+x\psi'(x)$ is increasing on $(0,\infty)$. In particular, for $x\geq \frac{1}{2}$, we have 
$$\psi(x)+x\psi'(x)\geq \psi\left(\frac{1}{2}\right)+\frac{1}{2}\psi'\left(\frac{1}{2}\right)=-\gamma-2\ln 2 +\frac{\pi^2}{4}\approx 0.503>0,$$ see relation (5.4.13) from Olver \textit{et al.} \cite{Olver},  
where $\gamma\approx 0.5772$ is the Euler--Mascheroni constant. 
In particular, one has that $\psi(rx)+rx\psi'(rx)\geq 0$ for every $  r\geq 1$ and $x\geq \frac{1}{2}.$

Let $g(x,r):=r\psi(rx)$. Then by the latter relation, one has that 
$$\frac{\partial}{\partial r}g(x,r)=\psi(rx)+rx\psi'(rx)\geq 0,\ \ \forall r\geq 1,\ x\geq \frac{1}{2}.$$ Therefore, $r\mapsto g(x,r)$ is increasing on $[1,\infty)$ for every $x\geq \frac{1}{2};$ accordingly, 
$$r\psi(rx)=g(x,r)\geq g(x,1)=\psi(x),\ \ \forall r\geq 1,\ x\geq \frac{1}{2},$$
 that is exactly the inequality we wanted to prove. 
%If $r\geq 1$ is fixed and  $h_r(x):=\frac{ \Gamma(rx)}{ \Gamma(x)}$, it follows by the latter property that 
%$$\frac{h_r'(x)}{h_r(x)}=r\psi(rx)-\psi(x)\geq 0,\ \ \forall x\geq \frac{1}{2},$$
%i.e. $h_r$ is increasing on $[\frac{1}{2},\infty).$
\hfill $\square$

	\begin{remark}\rm 
Although is not important for our arguments, we notice that the monotonicity property from  Proposition	\ref{Alzer-monoton-2} \textit{cannot} be extended to the whole interval  $(0,\infty)$. Having the previous proof, we believe that   $x\mapsto \frac{ \Gamma(rx)}{ \Gamma(x)}$ is increasing on the maximal interval $[x_0,\infty)$ for every $r\geq 1$ if and only if $x_0\approx 0.2161$ is the unique root of the equation $\psi(x)+x\psi'(x)=0.$ 
\end{remark}

\subsection{A parallelogram-type  inequality}

We conclude the paper by proving the following parallelogram-type  inequality that has been used in the proof of Theorem \ref{theorem-general}.

\begin{proposition}\label{elementary-inequality}
	If $q\in (1,2]$, then for every $v,w\in \mathbb R^k$ one has   
	\begin{equation}\label{elemi-1}
		|v+w|^q\leq (3-q)|v|^q+q|w|^{q-2}\langle v,w\rangle + |w|^q.
	\end{equation}
\end{proposition}

{\it Proof.} Let $|v|=r$, $|w|=R$ and $\langle v,w\rangle=rR\cos t$, $t\in [0,2\pi)$. Then \eqref{elemi-1} is equivalent to 
\begin{equation}\label{elemi-2}
	(r^2+R^2+2rR\cos t)^\frac{q}{2}\leq (3-q)r^q+qrR^{q-1} \cos t + R^q.
\end{equation}
If  either $r=0$ or $R=0$ or $q=2$, \eqref{elemi-2} trivially holds. In fact, in the case $q=2$ we have  the  parallelogram-rule making this inequality an identity. 

Thus, we  assume that $r,R>0$ and $q\in (1,2)$. 
The main tool used to prove  \eqref{elemi-2} is  Bernoulli's inequality
\begin{equation}\label{Bernoulli}
	(1+x)^\alpha\leq 1+\alpha x,\ \ \forall x\geq -1,\ \alpha\in [0,1].
\end{equation}
Since \eqref{elemi-2} is not symmetric in $r$ and $R$, we have to distinguish two cases.

\textit{Case 1}: $r\leq R.$  By applying \eqref{Bernoulli} for  
$ x :=\frac{r^2}{R^2}+2\frac{r}{R}\cos t \geq -1$  and  $ \alpha = \frac{q}{2} \in (0,1) ,$  it follows that 
\begin{eqnarray*}
R^q\left(1+\frac{r^2}{R^2}+2\frac{r}{R}\cos t\right)^\frac{q}{2}&\leq& R^q\left(1+\frac{q}{2}\left(\frac{r^2}{R^2}+2\frac{r}{R}\cos t\right)\right)\\&=&R^q+qrR^{q-1} \cos t+\frac{q}{2}R^{q-2}r^2.
\end{eqnarray*}
Since  $q< 2$ and $r\leq R$, it follows that 
$\frac{q}{2}R^{q-2}r^2\leq (3-q)r^q,$
 which ends the proof of \eqref{elemi-2}.

\textit{Case 2}: $r> R.$ Again by \eqref{Bernoulli},  it follows that 
\begin{eqnarray*}
	r^q\left(1+\frac{R^2}{r^2}+2\frac{R}{r}\cos t\right)^\frac{q}{2}&\leq& r^q\left(1+\frac{q}{2}\left(\frac{R^2}{r^2}+2\frac{R}{r}\cos t\right)\right)\\&=&r^q+qRr^{q-1} \cos t+\frac{q}{2}r^{q-2}R^2.
\end{eqnarray*}
Thus, in order to prove \eqref{elemi-2}, we have to show that 
\begin{equation}\label{elemi-3}
	(2-q)r^q+R^q-\frac{q}{2}R^2r^{q-2}+qrR(R^{q-2}-r^{q-2})\cos t\geq 0, \ \forall t\in [0,2\pi).
\end{equation}
Since $R<r$ and $q< 2$, we have that $R^{q-2}> r^{q-2}$; thus, instead of \eqref{elemi-3} it is enough to prove 
\begin{equation}\label{elemi-4}
	(2-q)r^q+R^q-\frac{q}{2}R^2r^{q-2}-qrR(R^{q-2}-r^{q-2})\geq 0.
\end{equation}
Let $s:=\frac{R}{r}\in (0,1)$ and consider the function $f:(0,1)\to \mathbb R$ defined by
$$f(s)=2-q+s^q-\frac{q}{2}s^2-qs(s^{q-2}-1).$$
To prove \eqref{elemi-4}, we are going to show that $f(s)\geq 0$ for every $s\in (0,1)$. Note that 
$0<2-q=f(0^+)<3-\frac{3q}{2}=f(1^-).$
 Since 
 $f'(s)=qs^{q-1}-qs-q(q-1)s^{q-2}+q,$ it follows that $f'(0^+)=-\infty$ and $f'(1^-)=q(2-q)>0.$ In particular, $f$ has a minimum point $s_0\in (0,1)$, thus  $f'(s_0)=0$, which is equivalent to 
 $s_0^q-s_0^2-(q-1)s_0^{q-1}+s_0=0.$
Thus, using the latter relation, it follows that for every $s\in (0,1)$ one has
\begin{eqnarray*}
f(s)&\geq& f(s_0)=2-q+s_0^q-\frac{q}{2}s_0^2-qs_0(s_0^{q-2}-1)\\&=&2-q+\left(1-\frac{q}{2}\right)s_0^2-s_0^{q-1}+(q-1)s_0\geq 0,
\end{eqnarray*}
where we used Young's inequality $s_0^{q-1}=1^{2-q}\cdot s_0^{q-1}\leq 2-q+(q-1)s_0$. \hfill $\square$
\vspace{0.5cm}

%			
%			
%			
%			
%			
%			\noindent {\bf Acknowledgments.} The authors thank S. Brendle, G. Carron,  L. Mazzieri, for stimulating conversations in the early stage of the manuscript. We also wish to thank E. Milman and S. Ohta for their comments and advice. \\
%			
%			
%			
%			
%			\noindent \textbf{Final statement.} The authors state that there is no conflict of interest.

	\noindent {\bf Acknowledgment.} 
	The authors thank the anonymous Referees for their  constructive comments on the manuscript.

		\end{document}